\documentclass[12pt]{article}
\usepackage{amssymb}
\usepackage{amsmath}
\usepackage{amsbsy}
\usepackage{amsthm}
\usepackage{rotating}
\font\BBb=msbm10 at 12pt
\renewcommand{\Bbb}[1]{\mbox{\BBb #1}}

\newcommand{\be}{\begin{equation}}
\newcommand{\ee}{\end{equation}}
\newcommand{\ba}{\begin{eqnarray}}
\newcommand{\ea}{\end{eqnarray}}
\newcommand{\ban}{\begin{eqnarray*}}
\newcommand{\ean}{\end{eqnarray*}}
\newcommand{\lrarrow}{-\!\!\!-\!\!\!-\!\!\!-\!\!\!-\!\!\!\rightarrow}
\newcommand{\pt}{\partial}

\newtheorem{theo}{Theorem}[section]
\newtheorem{prop}[theo]{Proposition}
\newtheorem{lem}[theo]{Lemma}
\newtheorem{cor}[theo]{Corollary}
\newtheorem{defi}[theo]{Definition}

\newtheorem{conj}{Conjecture}[section]
\begin{document}
\title{The coarse geometric Novikov conjecture and uniform convexity}
\author{Gennadi Kasparov and Guoliang Yu\footnote{Partially supported by NSF and NSFC.}}
\date{ }
\maketitle
\section{Introduction}

The classic Atiyah-Singer index theory of elliptic operators on
compact manifolds has been vastly generalized to higher index
theories of elliptic operators on noncompact spaces in the
framework of noncommutative geometry [5] by Connes-Moscovici for
covering spaces [8], Baum-Connes for spaces with  proper and
cocompact discrete group actions [2], Connes-Skandalis for
foliated manifolds [9], and  Roe for noncompact complete
Riemannian manifolds [33]. These higher index theories have
important applications to geometry and topology.

In the case of a noncompact complete Riemannian manifold, the
coarse geometric Novikov conjecture provides an algorithm to
determine when the higher index of an elliptic  operator on the
noncompact  complete Riemannian manifold is nonzero. The purpose
of this paper is to prove the coarse geometric Novikov conjecture
under a certain mild geometric condition suggested by Gromov [16].

Let $\Gamma$ be a metric space; let $X$ be a Banach  space. A map
$f : \Gamma \rightarrow X$ is said to be a (coarse) uniform
embedding [15] if there exist non-decreasing functions $\rho_1$
and $\rho_2$ from ${\Bbb R}_+ = [0,\infty )$ to ${\Bbb R}$ such
that
\begin{enumerate}
\item[(1)] $\rho_1 (d(x,y) ) \leq \Vert f(x) - f(y)\Vert \leq \rho_2
(d(x,y) )$ for all $x,y \in \Gamma$;
\item[(2)] $\lim_{r \rightarrow +\infty} \rho_i (r) = + \infty$ for
$i = 1,2$.
\end{enumerate}

In this paper, we prove the following result:
\begin{theo} Let $\Gamma$ be a discrete metric space with bounded geometry. If $\Gamma$ is
uniformly embeddable into a uniformly convex Banach space, then
the coarse geometric Novikov conjecture holds for $\Gamma$, i.e.,
the index map from
$\lim_{d\to\infty}K_*\left(P_d\left(\Gamma\right)\right)$ to
$K_*\left(C^* (\Gamma)\right)$ is injective, where $P_d(\Gamma)$
is the Rips complex of $\Gamma$ and  $C^* (\Gamma)$ is the Roe
algebra associated to  $\Gamma$.
 \end{theo}

  Recall that a discrete metric space $\Gamma$ is said
to have bounded geometry if $\forall ~r > 0$, $\exists ~N(r) > 0$
such that the number of elements in  $B(x,r)$ is at most $N(r)$
for all $x \in \Gamma$, where $B(x,r) = \{ y \in \Gamma : d(y,x)
\leq r\}.$ A Banach space $X$ is called uniformly convex if
$\forall\,\varepsilon
>0$, $\exists\, \delta > 0$ such that if $x,\,y\in S(X)$ and
$\|x-y\|\geq\varepsilon$, then $\left\|\frac{x+y}{2}\right\|<
1-\delta$, where $S(X)=\left\{x\in X,\,\|x\| =1\right\}$.

The coarse geometric Novikov conjecture implies that the higher
index of the Dirac operator on a uniformly contractible Riemannian
manifold is nonzero (recall that a Riemannian manifold is said to
be uniformly contractible if for every $r>0$, there exists $R\geq
r$ such that every ball with radius $r$ can be contracted to a
point in a ball with radius $R$). By Proposition 4.33 of [33] and
the Lichnerowicz argument, Theorem 1.1 implies the following
result:
\begin{cor}
Let $M$ be a Riemannian manifold with bounded geometry. If $M$
admits a uniform embedding into a uniformly convex Banach space
and is uniformly contractible, then $M$ cannot have uniformly
positive scalar curvature.
\end{cor}

In general, Gromov conjectures that a uniformly contractible
Riemannian manifold with bounded geometry cannot have uniformly
positive scalar curvature [16].

The possibility of using a uniform embedding into a uniformly
convex Banach space in order to study the Novikov conjecture was
suggested by Gromov [16]. The main new ideas in the proof of
Theorem 1.1 consist of a construction of a family of uniformly
almost flat Bott vector bundles over the uniformly convex Banach
space and a K-theoretic  finitization technique. The uniform
convexity condition is used in a crucial way to construct this
family of uniformly almost flat Bott vector bundles.

The coarse geometric Novikov conjecture is false if the bounded
geometry condition is removed   [39]. The coarse geometric Novikov
conjecture for bounded geometry spaces uniformly embeddable into
Hilbert space was proved in [40]. The proof of the Hilbert space
case makes the use of an algebra of the Hilbert space introduced
in [21] [22] [23].

 W. B. Johnson and N. L. Randrianarivony showed that $l_p$ $(p>2)$ does not
admit a (coarse) uniform embedding into a Hilbert space [27]. More
recently, M. Mendel and  A. Naor proved that $l_p$ does not admit
a (coarse) uniform embedding into $l_q$ if $p>q\geq 2$ [31].
 N. Brown and E. Guentner proved that every bounded
geometry space admits a (coarse) uniform embedding into a strictly
convex and reflexive Banach space [3] (recall that a Banach space
$X$ is said to be strictly convex if
$\left\|\frac{x+y}{2}\right\|<1$ for any two distinct unit vectors
$x$ and $y$ in $X$). It is an open question whether every
separable metric space admits a (coarse) uniform embedding into
some uniformly convex Banach space. N. Ozawa proved that expanders
don't admit a (coarse) uniform embedding into  a uniformly convex
Banach space with an unconditional basis [32]. We also would like
to mention the conjecture that if $M$ is a compact smooth
manifold, then any countable subgroup of the diffeomorphism group
$Diff(M)$ of $M$ admits a (coarse) uniform embedding into $C_p$
for some $p>1$, where $C_p$ is the Banach space of all
Schatten-$p$ class operators on a Hilbert space (recall that $C_p$
is uniformly convex for all $p>1$).

We remark that the K-theory for complex Banach algebras throughout
this paper is the 2-periodic complex topological K-theory.

This paper is organized as follows: In Section 2, we collect a few
facts about uniform convexity which will be used later in this
paper. In Section 3, we introduce (Banach) Clifford algebras over
a Banach space. In the case of a uniformly convex Banach space, we
use the (Banach) Clifford algebras  to construct a certain Bott
vector bundles over the Banach space and show that the Bott vector
bundles are uniformly almost flat in a certain Banach sense. As
suggested by Misha Gromov, Section 3 might be of independent
interest to experts in Banach space theory. In Section 4, we
briefly recall the coarse geometric Novikov conjecture and the
K-theoretic localization technique. Section 5, we introduce a
K-theoretic finitization technique. In Section 6, we use the Bott
vector bundles  to construct Bott maps in K-theory. Finally in
Section 7, we prove the main result of this paper.

In a separate paper,  we will show how uniform convexity can be
used to study K-theory for $C^*$-algebras associated to discrete
groups. In particular, we shall prove the Novikov conjecture for
groups uniformly embeddable into uniformly convex Banach spaces.

The authors would like to thank Misha Gromov, Erik Guentner, John
Roe and the referees for very helpful comments on earlier versions
of this paper which led to this revised version.

\section{Uniform convexity of Banach spaces}

In this section, we collect a few facts about uniformly convex
Banach spaces which will be used in this paper. A beautiful
account of the theory of uniformly convex Banach space can be
found in  Diestel's book [10].

For convenience of the readers, we give a proof of the following
classic result in the theory of Banach spaces.

\begin{prop} Let $X$ be a Banach space (over $\mathbb{R}$). Assume that $X^*$ is
uniformly convex.
\begin{enumerate}
\item[(1)] $\forall\,x \in S (X)$, $\exists$ a unique $x^* \in S(X^*)$ such that
$x^*(x)=1$.
\item[(2)] The map: $x \to x^*$, from $S(X)$ to $S(X^*)$ is uniformly continuous.
\end{enumerate}
\end{prop}
{\noindent {\it Proof:}} Let's first prove part (1) of Proposition
2.1 . Existence of $x^*$ follows from the Hahn-Banach theorem. Let $g$ be another element in $S(X^*) $
satisfying $g(x)=1$. We have $\frac{x^*+g}{2}(x)=1$. This implies
that $\|\frac{x^*+g}{2}\|=1$. Uniform convexity of $X^*$ implies
$g=x^*$.

Next we shall prove part (2) of the proposition. Given
$\varepsilon>0$, let $\delta>0$ be as in the definition of uniform
convexity of $X^*$. Assume that a pair of vectors $x$ and $y\in S(X)$
satisfies $\|x-y\|< \delta$. We have $\frac{x^*+y^*}{2}
(x)>1-\delta.$ This implies that $\| \frac{x^*+y^*}{2} \| >
1-\delta$. By uniform convexity of $X^*$, we have
$\|x^*-y^*\|<\varepsilon.$
\begin{flushright}
$\Box$
\end{flushright}

We remark that part (1) of Proposition 2.1 is still true under the
weaker condition that $X^*$ is strictly convex.

For any $x \in X$, we define $x^* \in X^*$ by:
\[
x^*= \left\{\begin{array}{ll} \|x\|\left(\frac{x}{\|x\|}\right)^*
& \mathrm{if}\ {\it x}\not= 0\,; \\ 0 &  \mathrm{otherwise}\,.
\end{array}\right.
\]

The following result is a consequence of Enflo's theorem [13]
and Asplund's averaging technique [1]. (See also [10], p. 87.) It  plays an
important role in the proof of the main result of this paper.

\begin{theo}
Let $\Gamma$ be a discrete metric space. If $\Gamma$ admits a
uniform embedding into a uniformly convex Banach space, then
$\Gamma$ is uniformly embeddable into a uniformly convex Banach
space $X$ such that its dual space $X^*$ is also uniformly convex.
\end{theo}

\section{(Banach) Clifford algebras and Bott vector bundles over Banach spaces}

In this section, we shall first introduce (Banach) Clifford
algebras over a Banach space. In the case of a uniformly convex
Banach space, we use the (Banach) Clifford algebras to construct
Bott vector bundles over the Banach space and show that the Bott
vector bundles are uniformly  almost flat in a certain Banach
sense.

Let $X$ be a Banach space over $\mathbb{R}$.  Let $V$ be a finite
dimensional subspace of $X$ and $V^*$ its dual space. We define a
pairing: $\left(V\oplus V^*\right)\times \left(V \oplus V^*\right)
\to \mathbb{R}$ by:
\[
q(x\oplus g),\,y\oplus h)=h(x)+ g(y)
\]
for all $ x\oplus g,\ y\oplus h \in V\oplus V^*$.

Let $W=V\oplus V ^*$ be given the norm:
 \[\|x\oplus g\|=\sqrt{\|x\|^2 +\|g\|^2}\]
 for all $x\oplus g\in W=V\oplus V ^*$.

Let $\otimes^n
W=\overbrace{W\otimes\cdots\otimes W}^n$ for $n\geq 1$ and
$\otimes^0 W=\mathbb{ R}$.

Endow $\otimes^0 W$ with the standard norm.
For  $n\geq 1$, endow $\otimes^n W$ with the following
norm:

$$ \| u\|= \sup \{ (\phi_1\otimes \cdots \otimes \phi_n )(u): \phi_k \in W^*, \|\phi_k\|\leq 1,\,\,\,\, 1\leq k\leq n\}$$
for all $u\in \otimes^n W$,
where $W^*$  is the dual (Banach) space of $W$.

Let
\[
T(W)=\left\{\underset{n=0}{\overset{\infty}{\oplus}} a_n\,:\,
a_n\in \otimes^n W\,,\, \sum_{n=0}^{\infty}\| a_n  \|< +\infty
\right\}
\]
be the tensor algebra.

Observe that $T(W)$ is a Banach algebra over
$\mathbb{R}$ with the following norm
$$\|\underset{n=1}{\overset{\infty}{\oplus}} a_n \|=\sum_{n=0}^\infty \|a_n \|.$$

Let $T_{\mathbb{C}}(W)$ be the complexification of the Banach
algebra $T(W)$. $T_{\mathbb{C}}(W)$ is a Banach algebra over
$\mathbb{C}$.

Let $I_{\mathbb{C}}(W)$ be the closed two-sided ideal in
$T_{\mathbb{C}}(W)$ generated by all elements of the form
\[
w_1\otimes w_2+ w_2\otimes w_1 + 2q(w_1,w_2),\quad w_1,w_2 \in W.
\]
The Clifford algebra $Cl(W)$ is defined as the quotient Banach
algebra:
\[
Cl(W)=T_{\mathbb{C}}(W)/I_{\mathbb{C}}(W).
\]

$Cl(W)$ is a finite-dimensional  complex Banach algebra with the
natural quotient norm.

Let $C_b(W,Cl(W))$ be the Banach algebra of all bounded continuous
functions on $W$ with values in $Cl(W)$, where the norm of each
element $f\in C_b(W,Cl(W))$ is defined by:
\[\|f\|=\sup_{w\in W} \|f(w)\|.\]

Let $C_0(W,Cl(W))$ be the Banach algebra of all continuous
functions on $W$ with values in $ Cl(W)$, vanishing at $\infty$,
where the norm on $C_0(W,Cl(W))$ is inherited from the norm on
$C_b(W,Cl(W))$.

Throughout the rest of the paper, let $X$ be a uniformly convex
Banach space over $\mathbb{R}$ such that its dual space $X^*$ is
also uniformly convex. We shall next construct a family of
uniformly almost flat representatives of the Bott generators in
the K-group $K_0(C_0(W,Cl(W)))$ for all finite dimensional
subspaces $V \subseteq X$.

We need to recall the concept of the index of a generalized
Fredholm operator in the context of Banach algebras.

 Let $B$ be a ($\mathbb{Z}_2$-)graded
unital complex Banach algebra and let $A$ be a
($\mathbb{Z}_2$-)graded ideal in $B$. Assume that the grading on
$B$ is induced by a grading operator $\varepsilon$ in $B$
satisfying $\varepsilon^2=1, \|\varepsilon\|=1$. We have
$B^{(0)}=\{b\in B: \varepsilon^{-1}b \varepsilon=b\}$ and
$B^{(1)}=\{b\in B: \varepsilon^{-1}b \varepsilon=-b\}$.

 Let $F\in B$ be an element of degree one
such that $F^2-1 \in A$.

Write $$\alpha=\left(\frac{1+\varepsilon}{2}\right) F
\left(\frac{1-\varepsilon}{2}\right),$$
$$\alpha'=\left(\frac{1-\varepsilon}{2} \right) F
\left(\frac{1+\varepsilon}{2}\right).$$

Let
\[
\omega= \left(\begin{array}{cc}1 & \alpha \\ 0 &
1\end{array}\right)\left(\begin{array}{cc}1 & 0 \\ -\alpha' &
1\end{array}\right)\left(\begin{array}{cc}1 & \alpha \\ 0 &
1\end{array}\right)\left(\begin{array}{cc}0 & -1 \\ 1 &
0\end{array}\right).
\]
We have
\[
\omega=\left(\begin{array}{cc}\alpha+(1-\alpha\alpha')\alpha &
-(1-\alpha\alpha')\\ 1-\alpha'\alpha & \alpha'\end{array}\right)
\]
\[
\omega^{-1}=\left(\begin{array}{cc}\alpha'& 1-\alpha'\alpha\\
-(1-\alpha\alpha') & \alpha+ \alpha
(1-\alpha'\alpha)\end{array}\right).
\]
Define
\[
index(F)=\omega\left(\begin{array}{cc} 1 & 0 \\ 0 & 0
\end{array}\right)\omega^{-1}\,.
\]

Clearly, $index(F)$ is an idempotent in $M_2(A)$.

 We have
\[
index(F)=\left(\begin{array}{cc}\alpha\alpha'+(1-\alpha\alpha')\alpha\alpha'
&
\alpha(1-\alpha'\alpha)+(1-\alpha\alpha')\alpha(1-\alpha'\alpha)\\
(1-\alpha'\alpha)\alpha' & (1-\alpha'\alpha)^2\end{array}\right).
\]
 We define
\[
Index(F)=[index(F)]-\left[\left(\begin{array}{cc} 1 & 0 \\ 0 & 0
\end{array}\right)\right]\in K_0(A)\,.
\]
where $K_0(A)$ is the $K$-group of $A$ considered as a Banach
algebra without grading.

Notice that $Index(F)$ is  an obstruction to the invertibility of $F$.

Let $\varphi$ be a continuous function on $\mathbb{R}$ such that
$0\leq \varphi(t)\leq 1$, $\exists\, 0<c_1<c_2$ satisfying
$\varphi(t)=0$ if $t\leq c_1$ and $\varphi(t)=1$ if $t\geq c_2$.

 Note that the uniform convexity of $X^*$ is equivalent to the
uniform smoothness of $X$ (see [10], p. 36), therefore, all linear
subspaces $V$ of $X$ have the property that both $V$ and $V^*$ are
uniformly convex. Let $V$ be a finite dimensional subspace of $X$.
For any $x\in V$, let $x^* \in X^*$ be defined as above. The
restriction of $x^*$ to $V$ is still denoted by $x^*$. Since $V$
is naturally isometric to $(V^*)^*$ (with its natural norm),
$(V^*)^*$ is uniformly convex. For any $h\in V^*$, we identify
$h^*\in (V^*)^*$ with an element (still denoted by $h^*$) in $V$.

Let $F_{V, \varphi},\ \in C_b(W,Cl(W))$ be defined by:
\[
F_{V,\varphi}(0\oplus 0)=0,\]
\[
F_{V,\varphi}(x\oplus h)=
\]
\[
\frac{\varphi\left(\sqrt{\|x\|^2+\|h\|^2}\right)}{
\sqrt{\|x\|^2+\|h\|^2+i(h(x)-x^*(h^*))}}\left(\frac{h^*\oplus x^*-x\oplus h}{2}
+i\frac{x\oplus h+h^*\oplus x^*}{2}\right),
\]
for all nonzero $x\oplus h \in W=V\oplus V^*$.

Note that  $\sqrt{\|x\|^2+\|h\|^2+i(h(x)-x^*(h^*))}$ is well
defined as a continuous complex-valued function of $x\oplus h$
since $h(x)-x^*(h^*)$ are real numbers.

Endow $  C_b (W, Cl (W))$ with the  grading induced by the natural
grading operator $\varepsilon$ of $Cl(W)$ (considered as a
constant function on $W$). It is not difficult to see that the norm of
the grading operator is $1$.

Clearly  $F_{V,\varphi}$ has degree one.  It is also
straightforward to verify that
\[
F^2_{V,\varphi}-1 \in C_0(W,Cl(W))\, ,
\]
where $1$ is the identity element of $C_b(W,Cl(W))$.

\begin{lem} \textrm{Index}($F_{V,\varphi}$) is a generator for $K_0(C_0(W,Cl(W)))$.
\end{lem}
{\noindent {\it Proof:}} Let $\|\cdot\|_0$ be the Euclidean norm
on $V$. We define a homotopy of norms on $V$ and $V^*$ by:
$$\|\cdot\|_t =\sqrt{t\|\cdot\|^2+ (1-t) \| \cdot\|_0^2}$$
for $t\in [0,1]$.

For each $x\in V$, let $x^{*,0}\in V^*$ be given by
$x^{*,0}(y)=<y,x>$ for all $y\in V$ , where $<\, , \,>$ is the
inner product on $V$ corresponding to the Euclidean structure.  We
define $x^{*,t}\in V^*$ by:
$$ x^{*,t} =tx^* +(1-t)x^{*,0}$$
for $t\in [0,1]$. For each $h\in V^*$, let $h^{*,0}\in V$ be given
by $h(x)=<x, h^{*,0}>$ for all $x\in V$.
 We define $h^{*,t} \in V$ by:
$$ h^{*,t}=t h^{*}+ (1-t) h^{*,0}$$
for $t\in [0,1]$.

We can define a homotopy $F_{V,\varphi}(t)$ by respectively
replacing $x^*$, $h^*$ and $\|\cdot\|$ with $x^{*,t}$, $h^{*,t}$,
and $\|\cdot\|_t$ in the definition of $F_{V,\varphi}$. We have
$$Index (F_{V,\varphi})= Index(F_{V,\varphi}(1))=Index
(F_{V,\varphi}(0))$$ in $K_0(C_0(W,Cl(W)))$. But it is
straightforward to verify that $Ind(F_{V,\varphi}(0))$ is the Bott
generator for $K_0(C_0(W,Cl(W)))$.
\begin{flushright}
$\Box$
\end{flushright}

We have the following proposition.
\begin{prop} Let
\[
index(F_{V,\varphi})=a_{V,\varphi}+\left(\begin{array}{cc} 1 & 0 \\
0 & 0
\end{array}\right)
\]
for some $a_{V,\varphi}\in M_2(C_0(W,Cl(W)))$.
\begin{enumerate}
\item[(1)] Given $\varphi$, there exists $R>0$ such that
$supp(a_{V,\varphi})\subseteq B_W(0,R)$ for any finite-dimensional
subspace $V\subseteq X$, where $supp(a_{V,\varphi})=\left\{\xi\in
W\,:\, a_{V,\varphi}(\xi)\not=0\right\}$ ($a_{V,\varphi}$ is
identified as an element
of $C_0(W,M_2(Cl(W)))$, and $B_W(0,R)$\\
$=\left\{ \xi\in W\,:\, \|\xi\|<R\right\}$,
($\|\xi\|=\sqrt{\|x\|^2+\|h\|^2}$\,\ if \,\ $\xi=x\oplus h \in W=V
\oplus V^*$).
\item[(2)] There exists $C>0$ such that $\|index(F_{V,\varphi})\|\leq C$ for
any finite dimensional subspace $V\subseteq X$,  where the norm of each
element $f\in C_b(W,Cl(W))$ is defined by:
$\|f\|=\sup_{w\in W} \|f(w)\|, $ and $M_2 (C_b (W, Cl(W)))= $ \linebreak
$C_b(W,Cl(W))\otimes M_2( \mathbb{C})$ is endowed with a Banach algebra tensor product norm
(for example, the projective tensor product norm).
\end{enumerate}
\end{prop}
{\noindent {\it Proof:}} We have the following identity:
$$\left(\frac{h^*\oplus x^*-x\oplus h}{2}
+i\frac{x\oplus h+h^*\oplus x^*}{2}\right)^2=\|x\|^2+\|h\|^2+i(h(x)-x^*(h^*))$$
for all $x\oplus h \in W=V \oplus V^*$.

Now part $(1)$  of Proposition 3.2 follows
from the definition of $index (F_{V, \varphi})$, the definition of
$\varphi$, and the above identity.

Recall that $V$ and $V^*$ are real Banach spaces. Hence
we have $$|  \sqrt{\|x\|^2+\|h\|^2+i(h(x)-x^*(h^*))}\,\,\,|\geq \sqrt{\|x\|^2+\|h\|^2}$$
for all $x\oplus h \in W=V \oplus V^*$.

This inequality, together the definitions of
$F_{V,\varphi}$ and $\varphi$ and the norm on $V\oplus V^*$, implies that
$$\| F_{V,\varphi}\|\leq 2 .$$

Now part (2) of Proposition 3.2 follows
from the above inequality, the definition of $index (F_{V, \varphi})$,  and the fact that the norm
of the grading operator in $C_b(W,Cl(W))$ is $1$.

\begin{flushright}
$\Box$
\end{flushright}

The concept of almost flat bundles has been successfully used to
study the Novikov conjecture and positive scalar curvature problem
in [6] [7] [17] [19] [26]. We shall introduce a slight variation
of almost flatness suitable for the purpose of this paper.

 Given a natural number $k$, real numbers $r>0,$ $\epsilon>0$ and a subspace $U$ of $W$, an
idempotent $p$ in $M_k(C_b (W, Cl (W)))$ is said to be
$(r,\epsilon)$-flat relative to $U$ if

$$ \left\| p(u_1)-p(u_2)\right\|<\epsilon$$
for any $u_1$ and $u_2$ in $W$ satisfying $u_1 -u_2\in U$ and
$\|u_1- u_2\|\leq r$, where $p$ is identified with an element in
$C_b (W, M_k(Cl(W))).$

 The following result says that the family of
idempotents $\left\{index(F_{V,\varphi})\right\}_V$ is uniformly
almost flat.

\begin{prop} $\forall\,r>0,\ \epsilon>0$, there exists $\varphi$
such that $index(F_{V,\varphi})$ is $(r,\epsilon)$-flat relative
to $V\oplus 0 \subseteq W$ for any finite dimensional subspace
$V\subseteq X$, where $ W=V\oplus V^*$.
\end{prop}

{\noindent {\it Proof:}} Given $N>0$, $\delta>0$, let $\varphi$ be
a smooth function on $\mathbb{R}$ satisfying $$ 0\leq \varphi (t)
\leq 1, \,\,\,\, |\varphi ' (t) |< \delta$$ for all $t\in
\mathbb{R}$,  $\varphi (t)=0 $ if $t\leq N$,  and there exists
$N'>N$ such that $\varphi (t)=1$ if $t>N'$.

By Proposition 2.1, given $\epsilon'>0$, there exists $N>0$
(independent of $V$) such that
$$ \left \| \frac{x_1\oplus h}{\sqrt{\|x_1\|^2 +\|h\|^2}}-\frac{x_2 \oplus h}{\sqrt{\|x_2\|^2
+\|h\|^2}}\right\|<\epsilon', $$
$$\left \| \frac{h^* \oplus x_1^*}{\sqrt{\|x_1\|^2 +\|h\|^2}}-\frac{h^*\oplus x_2^*}{\sqrt{\|x_2\|^2
+\|h\|^2}}\right\|<\epsilon' $$

for all $x_1, x_2 \in V$, and $h\in V^*$ satisfying
$\|x_1-x_2\|<r$, $\|(x_1, h)\|\geq N$.

Write
\[ x'=\frac{x}{\sqrt{\|x\|^2+\|h\|^2}},\]
\[ h'=\frac{h}{\sqrt{\|x\|^2+\|h\|^2}}\]
for all nonzero $x\oplus h \in W=V\oplus V^*$.

We have
\[
F_{V,\varphi}(x\oplus h)=
\]
\[
\frac{\varphi\left(\sqrt{\|x\|^2+\|h\|^2}\right)}{
\sqrt{1+i(h'(x')-(x')^*((h')^*))}}\left(\frac{((h')^*\oplus (x')^*)-x'\oplus h'}{2}
+i\frac{x'\oplus h'+((h')^*\oplus (x')^*)}{2}\right),
\]
for all nonzero $x\oplus h \in W=V\oplus V^*$.

By choosing $\epsilon'$  and $\delta$ small enough, Proposition
3.3 follows from the above facts,  part (2) of Proposition 2.1,
part (3) of Proposition 2.3, the definitions of $\varphi$,
$F_{V,\varphi}$, $index (F_{V,\varphi})$, and straightforward
estimates.
\begin{flushright}
$\Box$
\end{flushright}

\section{The coarse geometric Novikov conjecture and localization}

In this section, we shall briefly recall the coarse geometric
Novikov conjecture and the localization technique.

 Let $M$ be a locally compact metric space. Let $H_M$ be a separable
Hilbert space equipped with a faithful and non-degenerate
$*$-representation of $C_0 (M)$ whose range contains no nonzero
compact operator.

\begin{defi} {\rm ([33])} (1) The support of a bounded linear operator $T : H_M
\rightarrow H_M$, denoted by $supp(T)$,  is the complement of the
set of points $(m, m') \in M \times M$ for which there exist
$\phi$ and $\psi$ in $C_0 (M)$ such that
$$\psi T \phi = 0, \quad \phi (m) \not= 0 \quad \text{and} \quad \psi
(m') \not= 0;$$

(2) The propagation of a  bounded operator $T: H_M \rightarrow
H_M$, denoted by $propagation (T)$, is defined to be  $$\sup \{
d(m,m') : (m,m') \in \text{supp} \,(T)\} ;$$

(3) A bounded operator $T: H_M \rightarrow H_M$ is locally compact
if the operators $\phi T$ and $T\phi$ are compact for all $\phi
\in C_0 (M)$.
\end{defi}

\begin{defi} {\rm ([33]) } The Roe algebra $C^\ast (M)$ is the operator norm closure of
the $\ast$-algebra of all locally compact, finite propagation
operators acting on $H_M$.
\end{defi}

It should be pointed out that $C^\ast (M)$ is independent of the choice of $H_M$ up to a $\ast$-isomorphism
(cf. [25]).
Throughout the rest of this paper, we  choose $H_M$ in the definition of $C^\ast (M)$
to be $l^2(Z)\otimes H$, where $Z$ is a countable dense subset of
$M$, $H$ is a separable and infinite dimensional Hilbert space,
and $C_0(M) $ acts on $H_M$ by: $\phi (g\otimes h)= (\phi g)\otimes
h$ for all $\phi \in C_0(M), g\in l^2(Z), h\in H$ ( $\phi$ acts on
$l^2(Z)$ by pointwise multiplication ).

Let $\Gamma$ be a locally finite discrete metric space (a metric
space is called locally finite if every ball contains finitely
many elements).

\begin{defi} For each $d \geq 0$, the Rips complex $P_d (\Gamma )$ is the
simplicial polyhedron where the set of all vertices is $\Gamma$,
and a finite subset $\{ \gamma_0 ,\dots ,\gamma_n \} \subseteq
\Gamma$ spans a simplex iff $d (\gamma_i ,\gamma_j ) \leq d$ for
all $0 \leq i, j \leq n$.
\end{defi}

Endow $P_d (\Gamma )$ with the  spherical metric. Recall that the
spherical metric is the maximal metric whose restriction to each
simplex is the metric obtained by identifying the simplex with part
of a unit sphere endowed with the standard Riemannian metric.
The distance of a pair of points in different connected components
of $P_d(\Gamma)$ is defined to be infinity. The use of spherical metric is necessary
to avoid certain pathological phenomena
when $\Gamma$ does not have bounded geometry. If $\Gamma$ has bounded geometry, one
can instead use the Euclidean/simplicial metric on $P_d(\Gamma)$.

\begin{conj}  {\rm (The coarse geometric Novikov  conjecture)} If $\Gamma$ is a
discrete metric space with bounded geometry, then the index map
$Ind$ from $\lim_{d \rightarrow \infty} K_\ast (P_d (\Gamma ) )$
to $\lim_{d \rightarrow \infty} K_\ast (C^\ast (P_d (\Gamma )))$
is injective, where  $ K_\ast (P_d (\Gamma ) ) =KK_\ast (C_0 (P_d
( \Gamma), {\Bbb C})$ is the $K$-homology group of the locally compact space
$P_d(\Gamma)$.
\end{conj}

It should be pointed out that $\lim_{d\rightarrow \infty} K_\ast
(C^\ast (P_d (\Gamma )))$ is isomorphic to $ K_\ast (C^\ast
(\Gamma ))$.

Conjecture  is false if the bounded geometry condition is dropped
[39].   By Proposition 4.33 in  [33], the coarse geometric Novikov
conjecture implies Gromov's conjecture that  a uniformly
contractible Riemannian manifold with bounded geometry cannot have
uniformly positive scalar curvature  and the zero-in-the spectrum
conjecture stating that the Laplace operator acting on the space
of all $L^2$-forms of a uniformly contractible Riemannian manifold
has zero in its spectrum.

The localization algebra introduced in [38] will play an important
role in the proof of our main result.  For the convenience of the
readers, we shall briefly recall its definition and its relation
with K-homology.

\begin{defi} Let $M$ be a locally compact metric space.
 The localization algebra $C_L^\ast (M)$ is the norm-closure
of the algebra of all uniformly bounded and uniformly
norm-continuous functions $$a : [0,\infty ) \rightarrow C^\ast
(M)$$ satisfying
$$propagation \,(a(t)) \rightarrow 0$$
as $t \rightarrow \infty$.
\end{defi}

There exists a local index map ([38]):
$$\text{Ind}_L : K_\ast (M) \rightarrow K_\ast (C_L^\ast (M)) .$$

\begin{theo} {\rm ([38])} If $P$ is a locally compact and finite dimensional simplicial polyhedron
endowed with the spherical metric, then the local index map
$\text{Ind}_L: K_{\ast}(P) \rightarrow
K_{\ast}(C_{L}^{\ast}(P))$  is an isomorphism.
\end{theo}

For the convenience of readers, we give an overview of the proof
of the above theorem given in [38]. Given a locally compact and
finite dimensional simplicial polyhedron  $P$ endowed with
 the spherical metric,
let $P_{1}$ and $P_{2}$ be two simplicial sub-polyhedrons of $P$.
Endow $P_{1}$, $P_{2}$, $P_1\cup P_2$ and $P_1\cap P_2$
 with  metrics inherited from
$P$. The local nature of the localization algebra can be used to
prove a Mayer-Vietoris for the K-groups of the localization
algebras for $P_{1}\cup P_{2}$, $P_{1}$, $P_{2}$ and $P_{1}\cap
P_{2}$, and a certain (strong) Lipschitz homotopy invariance of the
K-theory of the localization algebra. Now Theorem 4.5 follows from
an induction argument on the dimension of skeletons of $P$ using
the Mayer-Vietoris sequence and the (strong) Lipschitz homotopy
invariance for K-theory of localization algebras.

The evaluation homomorphism $e$ from
 $C_L^\ast (M)$ to $C^\ast (M)$ is defined by:
$$ e(a)=a(0)$$
for all $a\in C_L^\ast (M)$.

In the definitions of $C^*(P_d(\Gamma))$ and $C^*_L
(P_d(\Gamma))$, we choose  a countable dense subset $\Gamma_d$ of
$P_d (\Gamma)$ in such a way that if $d'>d$, then $\Gamma_d
\subseteq \Gamma_{d'}$. Hence there are natural inclusion
homomorphisms from $C^*(P_d(\Gamma))$ to $C^*(P_{d'}(\Gamma))$ and
from $C^*_L(P_d(\Gamma))$ to $C^*_L (P_{d'}(\Gamma))$ when $d'>d$.

If $\Gamma$ is a locally finite discrete metric space, we have the
following commuting diagram:

\vspace{1.00cm}

$$\lim_{d \rightarrow \infty} K_\ast (P_d (\Gamma ))$$
$$\text{Ind}_L \,\swarrow \hskip 1.0in \searrow \,\text{Ind}$$
$$\lim_{d \rightarrow \infty} K_\ast (C_L^\ast (P_d (\Gamma )))
\overset{e_{\ast}}{\longrightarrow}\lim_{d \rightarrow \infty}
K_\ast (C^\ast (P_d (\Gamma ))).$$

Theorem 4.5
 implies that in order to prove the coarse geometric
Novikov conjecture, it is enough to show that
$$e_\ast : \lim_{d \rightarrow \infty} K_\ast (C_L^\ast (P_d (\Gamma )))
\rightarrow \lim_{d \rightarrow \infty} K_\ast (C^\ast (P_d
(\Gamma )))$$ is injective.

Let $P$ be a locally compact and finite dimensional simplicial
polyhedron endowed with the spherical metric. Let $Q$ be a
simplicial sub-polyhedron of $P$ with the metric inherited from
$P$.

Let $P\backslash Q=\{x\in P: x\notin Q\}.$  We define the relative K-homology group of $(P,Q)$ by:
$$K_* (P, Q)=K_*(P\backslash Q)=KK_*(
C_0(P\backslash Q), {\Bbb C}).$$

Endow $P\backslash Q$ with the metric inherited from $P$.
We define $C_L^*(\partial (P\backslash Q))$ to be the closed subalgebra of $C_L^*(P\backslash Q)$
generated by elements $a \in C_L^*(P\backslash Q)$ such that there exists
$c_t>0$  ($t\in [0,\infty)$) satisfying $\lim_{t\rightarrow
\infty} c_t=0$, and $supp (a(t))\subseteq \{(x,y)\in (P\backslash Q)\times (P\backslash Q):
d((x,y), Q\times Q)<c_t\}$ for all $t\in [0,\infty),$ where $c_t$
depends on $a$. Notice that $C_L^*(\partial (P\backslash Q))$ is  the closed two sided ideal of $C_L^*(P\backslash Q).$

Let $C^*_L (Q; P)$ be the closed subalgebra of $C_L^*(P)$
generated by all elements $a$ in $C_L^*(P)$ such that there exists
$c_t>0$  ($t\in [0,\infty)$) satisfying $\lim_{t\rightarrow
\infty} c_t=0$, and $supp (a(t))\subseteq \{(x,y)\in P\times P:
d((x,y), Q\times Q)<c_t\}$ for all $t\in [0,\infty),$ where $c_t$
depends on $a$. It is not difficult to see that $C^*_L (Q; P)$ is
a closed two-sided ideal of $C_L^*(P)$.

Observe that $C_L^*(P)/C^*_L (Q,P)$ is naturally isomorphic to
$C_L^*(P\backslash Q)/C_L^*(\partial(Q\backslash P))$.

Next we shall define  a  local index map

$$ Ind_L: \,\,K_*(P,Q)\rightarrow
K_*(C_L^*(P)/C_L^*(Q,P))\cong K_*(C_L^*(P\backslash Q)/C_L^*(\partial(Q\backslash P))).$$

Recall that if $M$ is a locally compact topological space, the $K$-homology groups
$K_{i}(M)=KK_{i}(C_{0}(M), \Bbb{C})$ $(i=0,1)$ are generated
by the following cycles modulo a certain equivalence relation [28]:
\\(1) a cycle for $K_{0}(M)$ is a pair $(H_{M}, F)$, where $H_{M}$
is a Hilbert space with a $*$-representation of $C_0(M)$ and $F$ is a bounded linear operator acting on $H_{M}$
such that $F^{\ast}F-I$ and $FF^{\ast}-I$ are locally compact,
and $\phi F -F \phi$ is compact for all $\phi \in C_{0}(M);$
\\(2) a cycle for $K_{1}(M)$ is a pair $(H_{M}, F)$, where $H_{M}$
is a Hilbert space with a $*$-representation of $C_0(M)$ and $F$ is a self-adjoint  operator acting on $H_{M}$
such that $F^{2}-I$ is locally compact,
and $\phi F -F \phi$ is compact for all $\phi \in C_{0}(M).$

Let $(H_{P\backslash Q}, F)$ be a cycle for $K_{0}(P;Q)=K_0(C_0 (P\backslash Q))$. Without loss of generality,
we can assume that
$H_{P\backslash Q}$ is a  non-degenerate $*$-representation of $C_0(P\backslash Q)$ whose range contains
no nonzero compact operator.

For each positive integer $n$, there exists a locally finite open cover
 $\{U_{n,i}\}_{i}$
 for $P\backslash Q$ such that
$diameter(U_{n,i})<1/n$ for all $i$.
Let $\{\phi_{n,i}\}_{i}$ be a continuous partition of unity subordinate
to $\{U_{n,i}\}_{i}$.
Define a family of operators $F(t)$ $(t\in [0, \infty))$
acting on $H_{ P\backslash Q}$ by:
$$ F(t)=  \sum_{i}( (n-t)
 \phi^{\frac{1}{2}}_{n,i} F \phi^{\frac{1}{2}}_{n,i}
+ (t-n+1) \phi^{\frac{1}{2}}_{n+1,i} F
\phi^{\frac{1}{2}}_{n+1,i})$$ for all positive integer $n$ and
$t\in [n-1,n]$, where the infinite sum converges in strong
topology. Lemma 2.6 of [39] implies that $F(t)$ is a bounded and
uniformly norm-continuous
 function from $[0,\infty)$ to the $C^{\ast}$-algebra of
 all bounded operators acting on $H_{P\backslash Q}$.
Notice that  $$propagation(F(t))\rightarrow 0~~ as
~~t\rightarrow \infty.$$
Using the above facts it is not difficult to see that
$F(t)$ is a multiplier of $C_{L}^{\ast}(P\backslash Q)/C_L^*(\partial(Q\backslash P))$ and $F(t)$ is a unitary
modulo $C^{\ast}_{L}(P\backslash Q)/C_L^*(\partial(Q\backslash P))$.
Hence $F(t)$ gives rise to an element $$[F(t)]\in K_0(C_L^*(P)/C_L^*(Q,P))\cong K_0(C_L^*(P\backslash Q)/C_L^*(\partial(Q\backslash P))).$$

We define the local index of the cycle $(H_{P\backslash Q}, F)$ to be
$[F(t)]$.

Similarly we can define the local index map

 $$Ind_L: \,\,K_1(P;Q) \rightarrow
K_1(C_L^*(P)/C_L^*(Q,P))\cong K_1(C_L^*(P\backslash Q)/C_L^*(\partial(Q\backslash P))).$$

\begin{prop}
$ Ind_L$ is an isomorphism from $K_* (P,Q)$ to
$K_*(C^*_L(P)/C_L^*(Q;P))$.
\end{prop}

{\noindent {\it Proof:}} We have the following commutative
diagram:

\vspace{1.00cm}

\begin{tabular}{ccccccc}
 $\rightarrow $&$K_*(Q)$ & $\rightarrow $ & $K_*(P)$ & $\rightarrow$ & $K_*(P, Q)$ &$\rightarrow $\\
  &$\downarrow$ &  & $\downarrow$ &  & $\downarrow $& \\
  $\rightarrow$ & $K_*(C^*_L(Q;P))$& $\rightarrow$ & $K_*(C_L^*(P))$ & $\rightarrow$&
  $K_*(C^*_L(P)/C_L^*(Q;P))$
  &$\rightarrow$

\end{tabular}

\vspace{1.00cm}

\begin{tabular}{cccccc}
   $\rightarrow$& $K_{*+1}(Q)$ & $\rightarrow$ & $K_{*+1}(Q)$ & $\rightarrow$&$\cdots$\\
    & $\downarrow$ &  & $\downarrow$ & & \\
    $\rightarrow$& $K_{*+1}(C_L^*(Q;P))$& $\rightarrow$ & $K_{*+1}(C_L^*(P)) $ & $\rightarrow$& $\cdots$,  \\
  \end{tabular}

\vspace{1.00cm}

 where the vertical map from $K_*(Q)$ to
$K_*(C^*_L(Q;P))$ is the composition of the local index map with
the  homomorphism from $K_*(C_L^*(Q))$ to
$K_*(C^*_L(Q;P))$ induced by the inclusion homomorphism from $C_L^*(Q)$ to $C^*_L(Q;P)$.

Now Proposition 4.6 follows from the five lemma, Theorem 4.5, and
the fact that the inclusion homomorphism from $K_*(C_L^*(Q))$ to
$K_*(C^*_L(Q;P))$ is an isomorphism (cf. Lemma 3.10 in [38]).

\begin{flushright}
$\Box$
\end{flushright}

\section{Finitization of K-theory }

In this section,  we introduce  algebras associated to a sequence
of finite metric spaces. The new algebras allow  us to localize
$K$-theory of the Roe algebra associated to an infinite metric
space to its finite metric subspaces (we call this process
finitization).

Let $\Gamma$ be a discrete metric space with bounded geometry, let
$\left\{F_n\right\}_{n=1}^\infty$ be a sequence of finite metric
subspaces of $\Gamma$.

We define
%\begin{eqnarray*}
$C^*_{alg}\left( \left\{F_n\right\}_n\right)$ to be the algebra
$$\bigg\{ \underset{n=1}{\overset{\infty}{\oplus}} a_n\,:\, a_n \in
C^*(F_n),  \  \left. \sup_{n}\|a_n\|<+\infty,\
\sup_n(propagation(a_n))<+\infty \right\}\,.$$
%\end{eqnarray*}

Endow $C^*_{alg}\left( \left\{F_n\right\}_n\right)$ with the
following norm:
\[
\left\|\underset{n=1}{\overset{\infty}{\oplus}} a_n\right\|=\sup_n\|a_n\|\,.
\]
\begin{defi} The $C^*$-algebra $C^*\left( \left\{F_n\right\}_n\right)$ is defined to be
the norm completion of $C^*_{alg}\left(
\left\{F_n\right\}_n\right)$.
\end{defi}

\begin{defi} The $C^*$-algebra $C^*(\left\{\pt_\Gamma F_n\right\}_n)$ is defined to be the
closed subalgebra of  $C^*\left( \left\{F_n\right\}_n\right)$
generated by elements
$\underset{n=1}{\overset{\infty}{\oplus}}a_n$ such that there
exists $r>0$ satisfying
\[
supp(a_n)\subseteq \left(F_n\cap
B_\Gamma(\Gamma-F_n,r)\right)\times\left(F_n\cap
B_\Gamma(\Gamma-F_n,r)\right)
\]
for all $n$, where $B_\Gamma(\Gamma-F_n,r)=\left\{x\in
\Gamma\,:\,d(x,\Gamma-F_n)<r\right\}$ if $\Gamma-F_n \neq
\emptyset$, and $B_\Gamma(\Gamma-F_n,r)=\emptyset$ if $\Gamma-F_n
= \emptyset$.
\end{defi}

It is easy to see that  $C^*(\left\{\pt_\Gamma F_n\right\}_n)$ is
a two-sided ideal in $C^*\left( \left\{F_n\right\}_n\right)$.

\begin{defi}
We define the $C^{\ast}$-algebra  $C^*(\left\{F_n,\pt_\Gamma
F_n\right\}_n)$ to be the quotient algebra  $C^*\left(
\left\{F_n\right\}_n\right)/C^*(\left\{\pt_\Gamma F_n\right\}_n)$.
\end{defi}

Throughout the rest of this  paper, let $\chi_n\in
l^\infty(\Gamma)$ be the characteristic function of $F_n$.

Define a homomorphism
\[
S\,:\, C^*(\left\{F_n,\pt_\Gamma F_n\right\}_n)\to
C^*(\left\{F_n,\pt_\Gamma F_n\right\}_n)
\]
by
\[
S[a_1\oplus a_2\oplus a_3 \oplus \ldots]=[\chi_1 a_2 \chi_1\oplus
\chi_2 a_3\chi_2\oplus  \chi_3 a_4 \chi_3\oplus \ldots]\,
\]
for all $[a_1\oplus a_2\oplus a_3\oplus \ldots]\in
C^*(\left\{F_n,\pt_\Gamma F_n\right\}_n).$

Denote by $\,C^*([0,1], \left\{F_n,\pt_\Gamma F_n\right\}_n)\,\,$
the $C^*$-algebra of all continuous functions on $[0,1]$ with
values in $C^*(\left\{F_n,\pt_\Gamma F_n\right\}_n).$

\begin{defi}
We define the $C^*$-algebra $C^*_S (\Gamma)$ to be $$ \left\{a\in
C^*([0,1], \left\{F_n,\pt_\Gamma F_n\right\}_n)\,:\, a(1) =S
(a(0)) \right\}.$$
\end{defi}

Define a $*$-homomorphism $j$: $C^*_S(\Gamma)\rightarrow
C^*(\left\{F_n,\pt_\Gamma F_n\right\}_n)\,$ by:

$$j(a)= a (0)$$
for every $a\in C^*_S(\Gamma).$

\begin{prop}
We have the following exact sequence:
%\begin{eqnarray*}
$$\cdots \rightarrow  K_{*+1}(C^*(\left\{F_n,\pt_\Gamma
F_n\right\}_n))\overset{(Id-S)_*}{\lrarrow}
K_{*+1}(C^*(\left\{F_n,\pt_\Gamma F_n\right\}_n)) \rightarrow
K_*(C^*_S(\Gamma))$$ $$ \overset{j_*}{\lrarrow}
K_*(C^*(\left\{F_n,\pt_\Gamma F_n\right\}_n))
\overset{(Id-S)_*}{\lrarrow} K_*(C^*(\left\{F_n,\pt_{\Gamma}
F_n\right\}_n)) \rightarrow \cdots\ ,$$ where $Id$ is the identity
homomorphism.
%\end{eqnarray*}
\end{prop}
{\noindent {\it Proof:}} $\,\,\,\,$ Let
$\,\,\,\,\,A_1=C^*([0,1],\left\{F_n,\pt_{\Gamma}
F_n\right\}_{n\,\, \,odd})\,\,\,\,$ be the algebra of all
continuous functions on the unit interval $[0,1]$ with values in
$C^*(\left\{F_n,\pt_{\Gamma} F_n\right\}_{n\,\, \,odd})$, let
$\,\,\,A_2= C^*([0,1], \left\{F_n,\pt_{\Gamma} F_n\right\}_{n\,\,
\,even})\,\,\,$ $\,\,\,\,$  be the algebra of all continuous
functions on the unit interval $[0,1]$ with values in
$C^*(\left\{F_n,\pt_{\Gamma} F_n\right\}_{n\,\, \,even})$, and
$B=C^*(\left\{F_n,\pt_{\Gamma} F_n\right\}_n)$.

Define $f_1:A_1\to B\,\,\,$ and $\,\,\, f_2:A_2\to B\,\,\,$ by:
\[
f_1[a_1\oplus a_3\oplus a_5\oplus a_7\oplus \ldots)]=[a_1
(1)\oplus \chi_2 a_3(0) \chi_2 \oplus a_3 (1) \oplus \chi_4 a_5
(0)\chi_4 \oplus a_5 (1)\oplus \ldots)]\,,
\]
\[
f_2[b_2 \oplus b_4\oplus b_6 \oplus \ldots)]=[\chi_1 b_2
(0)\chi_1\oplus b_2 (1)\oplus \chi_3 b_4 (0)\chi_3\oplus b_4
(1)\oplus \ldots)]\,.
\]
Let $P=\left\{(a,b)\,|\,f_1(a)=f_2(b)\right\}\subseteq A_1\oplus
A_2$. It is not difficult to see that $P$ is isomorphic to $C^*_S
(\Gamma)$. Notice that $f_1$ and $f_2$ are surjective
$*$-homomorphisms.  By the Mayer-Vietoris sequence [36], we have
the following exact sequence:
\begin{eqnarray*}
\cdots &\lrarrow & K_{*+1}(A_1)\oplus K_{*+1}(A_2)\lrarrow
K_{*+1}(B)\lrarrow K_{*}(P)\\
&\overset{(g_{1*},\,g_{2*})}{\lrarrow}& K_{*}(A_1)\oplus
K_{*}(A_2)\overset{f_{2*}-f_{1*}} \lrarrow K_{*}(B)\lrarrow
\cdots\ ,
\end{eqnarray*}
where $g_1(a,b)=a$, $g_2(a,b)=b$. This implies Proposition 5.5.
\begin{flushright}
$\Box$
\end{flushright}

We   define a map
\[
\chi\,:\,C^*(\Gamma)\to C^*_S(\Gamma)
\]
by
\[
\chi(a)=[\underset{n=1}{\overset{\infty}{\oplus}}\chi_n a \chi_n
]\,
\]
for all $a\in C^*(\Gamma),$ where $\chi(a)$ is viewed as a
constant function on $[0,1]$.

It is not difficult to prove the following lemma.

\begin{lem}
$\chi$ is a $*$-homomorphism.
\end{lem}

Proposition 5.5 allows us to localize $\chi_* (K_* (C^*(\Gamma)))$
to finite subspaces of $\Gamma$. It is not known if $\chi_*$ is an
isomorphism.

Next we shall introduce twisted versions of the  algebras of
definitions 5.1, 5.2 and 5.3.

Let $f:\Gamma\to X$ be a uniform embedding. For each $n\in
\mathbb{N}$ (the set of all natural numbers),  let $V_n\subseteq
X$ be a finite-dimensional subspace such that $f(F_n)\subseteq
V_n$. Let $W_n=V_n\oplus V_n^*$.

Let $H$ be an infinite dimensional
separable Hilbert space, and let $H_{F_n}=l^2(F_n)\otimes H$ be the Hilbert space as in the definition of $C^*(F_n)$.

Recall that, for each $m\geq 1$, we define a norm on $\otimes ^m W_n$ by:

$$\| u\|=\sup \{ (\phi_1 \otimes \cdots \phi_m) (u): \,  \phi_k \in W_n^*,\,\,
\|\phi_k\|\leq 1, \,\, 1\leq k\leq m\}$$
for all $u \in \otimes ^m W_n$,
where
$W_n^*$ is the dual (Banach) space of $W_n$.

We define $T(W_n)$ to be the Banach space

$$\{ \underset{m=0}{\overset{\infty}{\oplus}} u_m: u_m \in \otimes ^m W_n, \sum _{m=0}^{\infty} \|u_m\|<\infty\}$$
endowed with the norm

$$\| \underset{m=0}{\overset{\infty}{\oplus}}  u_m\| = \sum _{m=0}^{\infty} \|u_m\|.$$

Let $T_{ \mathbb{C}}(W_n)$ be the complexification of $T(W_n)$.

For each  $m\geq 1$, we define a norm on $H_{F_n}\otimes (\otimes ^m W_n)$ by:

$$\| u\|=$$
$$\sup \{ (\phi_0 \otimes \phi_1 \otimes \cdots \phi_m)(u): \, \phi_0\in (H_{F_n})^*, \| \phi_0 \|=1,\, \phi_k \in W_n^*,\,\,
\|\phi_k\|\leq 1 ,\,\, 1\leq k\leq m\}$$
for all $u \in H_{F_n}\otimes (\otimes ^m W_n)$,
where $(H_{F_n})^*$ is the dual (Hilbert) space of $H_{F_n}$, and
$W_n^*$ is the dual (Banach) space of $W_n$.

We define $H_{F_n}\otimes T(W_n)$ to be the Banach space

$$\{ \underset{m=0}{\overset{\infty}{\oplus}} u_m: u_m \in H_{F_n}\otimes (\otimes ^m W_n), \sum _{m=0}^{\infty} \|u_m\|<\infty\}$$
endowed with the norm

$$\| \underset{m=0}{\overset{\infty}{\oplus}} u_m\| = \sum _{m=0}^{\infty} \|u_m\|.$$

Let $H_{F_n} \otimes T_{\mathbb{C}} (W_n)$ be the complexification of the Banach space
$H_{F_n}\otimes T(W_n)$.

Let $H_{F_n} \otimes I_{\mathbb{C}} (W_n)$ be the closed (complex) Banach
subspace of $H_{F_n} \otimes T_{\mathbb{C}} (W_n)$ spanned by all elements of the form

$$ h \otimes (v_1\otimes ((w_1\otimes w_2) +(w_2\otimes w_1 )+ 2q(w_1,w_2))\otimes v_2),$$
where $h\in H_{F_n}$, $w_1\in W_n$, $w_2\in W_n$,
 $v_1\in \otimes ^m W_n$ for some $m$,  $v_2\in \otimes ^l W_n$ for some $l$,  $q$ is the quadratic form
 on $W_n$ defined in Section 3, and $ I_{\mathbb{C}} (W_n)$ is as  in Section 3.

Now we define the Banach space $H_{F_n}\otimes Cl(W_n)$ to be quotient Banach space:

$$H_{F_n}\otimes Cl(W_n)= (H_{F_n}\otimes T_{\mathbb{C}} (W_n))/(H_{F_n} \otimes I_{\mathbb{C}} (W_n)),$$
where $Cl (W_n)$ is as in Section 3.

Let $C_0(W_n, H_{F_n}\otimes Cl(W_n))$ be the Banach space of all continuous functions
on $W_n$ with values in $H_{F_n}\otimes Cl(W_n)$, vanishing  at $\infty$,
 where the norm of each
element $\xi\in C_0(W_n, H_{F_n}\otimes Cl(W_n)) $ is defined by:
\[\|\xi \|=\sup_{w\in W_n} \|\xi (w)\|.\]

Let $C^* (F_n)\otimes_{alg} C_0 (W_n, Cl(W_n))$ be the algebraic tensor product of $C^*(F_n)$
with $C_0(W_n,  Cl(W_n))$. We shall construct a representation of  the algebra $C^* (F_n)\otimes_{alg} C_0 (W_n, Cl(W_n))$
on the Banach space $C_0(W_n, H_{F_n}\otimes Cl(W_n))$  as follows.

For each $T\in C^* (F_n)$, let   $T\otimes 1 $ be the bounded operator on  $H_{F_n}\otimes Cl(W_n)$
defined by:

$$(T\otimes 1)(h\otimes \eta)= (Th)\otimes \eta$$
for all $h\in H_{F_n}$ and $ \eta \in Cl(W_n)$.

We define a bounded operator $\psi (T)$ on
$C_0(W_n, H_{F_n}\otimes Cl(W_n))$ by:

$$ ((\psi (T))(\xi))(w) = (T\otimes 1) (\xi (w))$$
for all $\xi \in  C_0(W_n, H_{F_n}\otimes Cl(W_n))$ and $w\in W_n$.

For every element $v\in Cl(W_n)$, let $1\otimes v$ be the bounded operator
on $H_{F_n}\otimes Cl(W_n)$
defined by:

$$(1\otimes v)(h\otimes \eta)= h\otimes (v\eta)$$
for all $h\in H_{F_n}$ and $ \eta \in Cl(W_n)$, where $v\eta$ is the product of $v$ with $\eta$
in $Cl(W_n)$.

For every $g\in C_0 (W_n, Cl(W_n))$, we define a bounded operator $M_g$ on
 $C_0(W_n, H_{F_n}\otimes Cl(W_n))$ by:

 $$(M_g (\xi)) (w) = (1\otimes g(w)) (\xi (w)) $$
 for all $\xi \in  C_0(W_n, H_{F_n}\otimes Cl(W_n))$ and $w\in W_n$.

 We now construct a representation of  the algebra $C^* (F_n)\otimes_{alg} C_0 (W_n, Cl(W_n))$
on the Banach space $C_0(W_n, H_{F_n}\otimes Cl(W_n))$ by:

$$(T\otimes g) \xi = \psi(T) (M_g \xi) $$
for all $T\in C^* (F_n)$, $g\in C_0 (W_n, Cl(W_n))$, and   $\xi \in  C_0(W_n, H_{F_n}\otimes Cl(W_n))$.

By the definitions of the operators $\psi(T)$, $M_g$ and the Banach space norm on $C_0(W_n, H_{F_n}\otimes Cl(W_n))$,
it is not difficult to verify
 $$\| T\otimes g\|\leq \|T\| \,\,\|g\|.$$

We now define the Banach algebra $C^*(F_n,V_n)$ to be the closure of algebraic
tensor product of $C^*(F_n)\otimes_{alg} C_0 (W_n, Cl(W_n))$ under the operator norm given  by the above
representation of the algebra $C^* (F_n)\otimes_{alg} C_0 (W_n, Cl(W_n))$
on the Banach space $C_0(W_n, H_{F_n}\otimes Cl(W_n))$.

We remark that
$C^*(F_n,V_n)$ is carefully defined    for the purpose of constructing the Bott map later on
(cf. the definition of the Bott map before Lemma 6.1).

We define the Banach space $C_0 (W_n, H\otimes Cl (W_n))$ in a way similar to  the  above definition
of the Banach space $C_0(W_n, H_{F_n}\otimes Cl(W_n))$.
Let $K$ be the $C^*$-algebra of all compact operators on the Hilbert space $H$,  let $K\otimes_{alg} C_0(W_n, Cl(W_n))$
be the  algebraic tensor product of $K$ with $C_0(W_n, Cl(W_n))$.
We can construct a
representation of the algebra
$K\otimes_{alg} C_0(W_n, Cl(W_n))$ on the Banach space $C_0 (W_n, H\otimes Cl (W_n))$
in way similar to the construction of the above representation of the algebra $C^* (F_n)\otimes_{alg} C_0 (W_n, Cl(W_n))$
on the Banach space $C_0(W_n, H_{F_n}\otimes Cl(W_n))$.

We define
$K\otimes C_0(W_n, Cl(W_n))$ to be operator norm closure of the algebra $K\otimes_{alg} C_0(W_n, Cl(W_n))$,
where the operator norm is given by a representation of the algebra $K\otimes_{alg} C_0(W_n, Cl(W_n))$
on the Banach space $C_0 (W_n, H\otimes Cl (W_n)).$

We can identify  $C^*(F_n, V_n)$ with the algebra of all functions
on $F_n\times F_n$ with values in $K\otimes C_0(W_n, Cl(W_n))$
with the following convolution product:
$$ (a\cdot b)(x,y) =\sum_{z\in F_n} a(x,z)\, b(z,y)$$
for all $a$ and $b$ in $C^*(F_n, V_n)$, and   $(x,y)\in F_n \times
F_n$.

For any $a \in C^*(F_n, V_n)$, we define
$$ supp(a)=\{ (x,y)\in F_n\times F_n \,:\, a(x,y)\neq 0\},$$
$$propagation(a)=sup\{d(x,y): (x,y)\in supp(a)\}.$$

For any $a\in C^*(F_n, V_n)$ and $(x,y)\in F_n \times F_n$, we can
identify $a(x,y)$ as a function from $W_n$ to $K\otimes Cl(W_n)$
and define
$$support (a(x,y))= \{\xi\in W_n \,:\, (a(x,y))(\xi)\neq 0\},$$
where $ K\otimes Cl(W_n)$ is the operator norm closure of $K\otimes_{alg} Cl(W_n)$
( the operator norm is given by the natural representation of $K\otimes_{alg} Cl(W_n)$
on the Banach space $H\otimes Cl(W_n)$).

Let

{\setlength\arraycolsep{2pt}\begin{eqnarray*} C^*_{alg}\left(
\left\{F_n,V_n\right\}_n\right)&=&\bigg\{\underset{n=1}{\overset{\infty}{\oplus}}
a_n \,:\,a_n\in C^*(F_n,V_n),\ \sup_n\|a_n\|<+\infty,\ \exists \,
r>0\ \\ && \ \ \ \mathrm{such\ that}\ propagation(a_n)<r\ \
\mathrm{for\ all}\, n,\ \mathrm{and}\ \exists \, R>0\\ && \ \ \
\mathrm{such\ that}\ support(a_n(x,y))\subseteq B_{W_n}(f(x)\oplus
0,R)\ \ \mathrm{for\ all}\, n\bigg\}.
\end{eqnarray*}}

Endow $C^*_{alg}(\{F_n,V_n\}_n)$ with the norm
\[
\left\|\underset{n=1}{\overset{\infty}{\oplus}} a_n\right\|=\sup_n\|a_n\|.
\]
\begin{defi}
We define the Banach algebra $C^*(\left\{F_n,V_n\right\}_n)$ to be
the norm completion of $C^*_{alg}(\left\{F_n,V_n\right\}_n)$.
\end{defi}
\begin{defi}
$C^*(\left\{\partial_\Gamma F_n,V_n\right\}_n)$ is defined to be
the closed subalgebra of $C^*(\left\{F_n,V_n\right\}_n)$ generated
by elements $\oplus_{n=1}^\infty a_n$ such that $\exists\; r>0$
satisfying
\[
supp(a_n)\subseteq (F_n\cap B_\Gamma(\Gamma-F_n,r))\times(F_n\cap
B_\Gamma( \Gamma-F_n,r))\,.
\]
\end{defi}
Note that $C^*(\left\{\partial_\Gamma F_n,V_n\right\}_n)$ is a
closed two sided ideal of $C^*(\left\{F_n,V_n\right\}_n)$.
\begin{defi}
We define $C^*(\left\{F_n,\partial_\Gamma F_n ,V_n\right\}_n)$ to
be the quotient algebra
\[
C^*(\left\{F_n,V_n\right\}_n)/C^*(\left\{\partial_\Gamma
F_n,V_n\right\}_n).
\]

\end{defi}

\section{The Bott maps}

In this section, we use the family of uniformly almost flat vector
bundles introduced in Section 3 to construct certain Bott maps.
These Bott maps play a crucial role in the proof of the main
result of this paper.

We shall first describe a difference construction in $K$-theory of
Banach algebras. Let $B$ be a unital Banach algebra, let $A$ be a
closed two sided ideal in $B$. Let $p$ and $q$ be idempotents in
$B$ such that $p-q\in A$. We shall define a difference element
$D(p,q)\in K_0(A)$ associated to the pair $p$ and $q$. When $A$ and $B$ are
$C^*$-algebras, the difference construction described here is compatible
with the difference construction in  KK-theory (as explained below).

Let
\[
Z(p,q)=\left(\begin{array}{cccc}q&0&1-q&0 \\
1-q&0&0&q\\0&0&q&1-q\\0&1&0&0\end{array}\right)\,.
\]
We have
\[
(Z(p,q))^{-1}=\left(\begin{array}{cccc}q&1-q&0&0 \\
0&0&0&1\\1-q&0&q&0\\0&q&1-q&0\end{array}\right)\,.
\]
Define
\[
D_0 (p,q)=(Z(p,q))^{-1}\left(\begin{array}{cccc}p&0&0&0 \\
0&1-q&0&0\\0&0&0&0\\0&0&0&0\end{array}\right)Z(p,q)\,.
\]
Let
\[
p_1=\left(\begin{array}{cccc}1&0&0&0 \\
0&0&0&0\\0&0&0&0\\0&0&0&0\end{array}\right)\,.
\]

Notice that $D_0(p,q)$ is an element in $M_4 (A^+)$ and  $D_0(p,q)=p_1$ modulo $M_4(A)$.

We define
\[
D(p,q)=[D_0(p,q)]-[p_1]
\]
in $K_0(A)$.

Next we shall explain that, when $A$ and $B$ are
$C^*$-algebras, the difference construction described above is compatible
with the difference construction in  KK-theory.

Recall that the difference element in $KK( \mathbb{C}, A)\cong K_0(A)$ is represented by the KK module
$(E, \phi, F)$, where $E=A\oplus A$ is the Hilbert module over $A$ with the inner product:
$$< (a_0\oplus a_1), (b_0\oplus b_1)>= a_0^*b_0 \oplus a_1^*b_1$$ for all $a_0\oplus a_1$ and
$b_0\oplus b_1$ in $A\oplus A$, $\phi=\phi_0\oplus \phi_1$ is the homomorphism from $\mathbb{C}$ to
$B(E)$ defined by: $$ (\phi_0 (c))a= cpa,$$
$$(\phi_1(c)) a= cqa$$ for all $c\in \mathbb{C}$ and $a\in A$, and $F$ is the operator acting on $E$
defined by:
\[
F=\left(\begin{array}{cc}0&1\\
1&0\end{array}\right)\,.
\]

Let $E'_0=E'_1= A\oplus A\oplus A \oplus A$ be the Hilbert module over $A$ with the inner product:
$$<a_1\oplus a_2\oplus a_3\oplus a_4, b_1\oplus b_2\oplus b_3\oplus b_4>= a_1^*b_1\oplus a_2^*b_2\oplus a_3^*b_3\oplus a_4^*b_4$$
for all $a_1\oplus a_2\oplus a_3\oplus a_4$ and $ b_1\oplus b_2\oplus b_3\oplus b_4$ in $A\oplus A\oplus A \oplus A$.

Let $E'=E'_0\oplus E'_1$ and $F'$ be the operator acting on $E'$ defined by:
\[
F'=\left(\begin{array}{cc}0&I\\
I&0\end{array}\right)\,,
\]
where $I$ is the identity element in $M_4(\mathbb{C})$.

Let $\phi'_0$ and $\phi'_1$ be respectively homomorphisms from $\mathbb{C}$ to $B(E'_0)$  and
$B(E'_1)$ defined by:
\[ (\phi'_0 (c)) v_0=
c\left(\begin{array}{cccc}p&0&0&0 \\
0&1-q&0&0\\0&0&0&0\\0&0&0&0\end{array}\right) v_0
\]
for all $c\in \mathbb{C}$ and $v_0\in E'_0$, and
\[ (\phi'_1 (c)) v_1=
c\left(\begin{array}{cccc}q&0&0&0 \\
0&1-q&0&0\\0&0&0&0\\0&0&0&0\end{array}\right) v_1
\]
for all $c\in \mathbb{C}$ and $v_1\in E'_1$.

Let $\phi'=\phi'_0\oplus \phi'_1$.
It is not difficult to see that
$(E, \phi, F)$ is equivalent to $(E', \phi', F')$ as KK modules in
$KK(\mathbb{C}, A)$. This, together with identity $(Z(p,q)^{-1}\oplus Z(p,q)^{-1})F'(Z(p,q)\oplus Z(p,q))=F'$,
implies  that
$(E, \phi, F)$ is equivalent to $(E', (Z(p,q))^{-1}\oplus (Z(p,q))^{-1}) \phi' (Z(p,q)\oplus Z(p,q)), F')$ as Kasparov modules in
$KK(\mathbb{C}, A)$, where $Z(p,q)$ is defined as above.

By the formula of $(Z(p,q))^{-1}$ as described above, we can verify
\[p_1=Z(p,q))^{-1}\left(\begin{array}{cccc}q&0&0&0 \\
0&1-q&0&0\\0&0&0&0\\0&0&0&0\end{array}\right)Z(p,q)\, ,\]
where $p_1$ is defined as above.

It follows that $(E,\phi, F)$ is equivalent to $D(p,q)$ in $KK(\mathbb{C}, A)\cong K_0(A)$.

Let $C^* (F_n)\otimes_{alg} C_0(W_n, Cl(W_n))^+$ be the
algebraic tensor product of $C^* (F_n)$ with $C_0 (W_n, Cl(W_n))^+$,
where $C_0(W_n, Cl(W_n))^+$ is obtained from $C_0(W_n, Cl(W_n))$ by
adjoining the identity.

We can construct a representation of the algebra $C^* (F_n)\otimes_{alg} C_0(W_n, Cl(W_n))^+$
on the Banach space $C_0 (W_n, H_{F_n}\otimes Cl (W_n))$ exactly in the same way as the representation
of  the algebra $C^* (F_n)\otimes_{alg} C_0(W_n, Cl(W_n))$ on the Banach space $C_0 (W_n, H_{F_n}\otimes Cl (W_n))$ .
 We define   the Banach algebra $C^*(F_n, V_n^+)$ to be the  closure of
the  algebra $C^* (F_n)\otimes_{alg} C_0(W_n, Cl(W_n))$ under the operator norm given
by the representation of the algebra $C^* (F_n)\otimes_{alg} C_0(W_n, Cl(W_n))^+$
on the Banach space $C_0 (W_n, H_{F_n}\otimes Cl (W_n)).$

We define
\begin{eqnarray*}
C^*_{alg}(\left\{F_n,V_n^+\right\}_n)&=&\bigg\{\underset{n=1}{\overset{\infty}{\oplus}}
a_n\,:\, a_n\in C^*(F_n, V_n^+),\,\sup_n \|a_n\|<+\infty, \\ && \
\ \ \exists\;  r>0\ \mathrm{such\ that}\ propagation(a_n)<r\
\mathrm{for \ all}\,\, n, \\ && \ \ \ \exists\ R>0\ \mathrm{such\
that}\ a_n(x,y)=c_n I+ b_n(x,y) \\ &&  \ \ \ \mathrm{for\
some}\ c_n \in  \mathbb{C},\, b_n \in  C^*(F_n, V_n) \\
&& \ \ \ \mathrm{satisfying}\ support(b_n(x,y))\subseteq
B_{Wn}(f(x)\oplus 0,R)\bigg\}\, ,
\end{eqnarray*}
where $I$ is the identity element in $C_0(W_n, Cl (W_n))^+$.

Endow $C^*_{alg}(\left\{F_n,V_n^+\right\}_n)$ with the norm:
$$\|\underset{n=1}{\overset{\infty}{\oplus}} a_n \|=\sup_n\|a_n\|.$$

Let $C^*(\left\{F_n,V_n^+\right\}_n)$ be the norm completion of
$C^*_{alg}(\left\{F_n,V_n^+\right\}_n)$. Notice that $C^*(\left\{F_n,V_n\right\}_n)$
is a closed two sided ideal of
$C^*(\left\{F_n,V_n^+\right\}_n)$.

We can similarly define
$C^*(\left\{\partial_\Gamma F_n,V_n^+\right\}_n)$ and
$C^*(\left\{F_n,
\partial_\Gamma F_n,V_n^+\right\}_n)$.

Next we shall define a Bott map
\[
\beta_{\left\{V_n\right\}}:K_0(C^*(\left\{F_n\right\}_n))\to
K_0(C^*(\left\{F_n,V_n\right\}_n))\,.
\]

Let $e=\underset{n=1}{\overset{\infty}{\oplus}} e_n\in
C^*(\left\{F_n\right\}_n)$ be an idempotent representing an
element in $K_0(C^*(\left\{F_n\right\}_n))$. Given any $\delta>0$,
$\exists\;e'=\underset{n=1}{\overset{\infty}{\oplus}} e_n'\in
C^*_{alg}(\left\{F_n\right\}_n) $ such that $\|e'-e\|<\delta$.

Define
\[
 p_0(e')=\underset{n=1}{\overset{\infty}{\oplus}} a_n \in
 C^*_{alg}(\left\{F_n,V_n^+\right\})\otimes M_2 (\mathbb{C})
\]
by:
\[
(a_n(x,y))(v\oplus h)=e'_n(x,y)\otimes(index
(F_{V_n,\varphi}))((v+f(x))\oplus h)
\]
for all $x,y \in F_n$ and $v\oplus h \in W_n=V_n\oplus V^*_n$,
where $index(F_{V_n, \varphi})$ is defined as in Proposition 3.3. Notice that
$p_0 (e')$ depends on $\varphi$.

The finiteness of  $\sup_{n} \| a_n \|$  follows from part (1) of
Lemma 6.1 below.  This, together with Proposition 3.2, implies that $ p_0 (e')$ is an element
in $C^*_{alg}(\left\{F_n,V_n^+\right\})\otimes M_2 (\mathbb{C}).$

\begin{lem}
\begin{enumerate}
\item[(1)] Let $b=\underset{n=1}{\overset{\infty}{\oplus}} b_n\in
C^*_{alg}(\left\{F_n\right\}_n),$ $\phi \in C_0(W_n, Cl (W_n))^+ $.
 Let  $c_n(b,\phi)=  \underset{n=1}{\overset{\infty}{\oplus}} c_n \in
C^*_{alg}(\left\{F_n, V_n^+\right\}_n)$  be defined  by:
$$ (c_n(x,y))(v\oplus h)=b_n(x,y)\otimes \phi ((v+f(x))\oplus h)$$
for all $x,y \in F_n$ and $v\oplus h\in W_n=V_n\oplus V_n^*$, where $\phi$ is identified
with a function on $W_n$ with values in $Cl(W_n)$.  We have

$$\| c(b, \phi)\|\leq \|b\|  \,\, \|\phi\|.$$

\item[(2)] Let $e,\,\, e'$ and $\delta$ be as above.
If $\Gamma$ has bounded geometry, then $\forall
\epsilon>0,\exists\ \delta>0$ and $\varphi$ such that
$$\|(p_0(e'))^2-p_0(e')\|<
\epsilon,$$ where $\varphi$ is as in the definition of $F_{V_n,
\varphi}$, both $\varphi$ and $\delta$ are independent of $n$.
\end{enumerate}
\end{lem}
{\noindent {\it Proof:}}

For each $m\geq 1 $, recall that the norm on $\otimes ^m W_n$ is defined by:

$$\| u\|=\sup \{ (\lambda_1 \otimes \cdots \lambda_m) (u): \,  \lambda_k \in W_n^*,\,\,
\|\lambda_k\|\leq 1, \,\, 1\leq k\leq m\}$$
for all $u \in \otimes ^m W_n$,
where
$W_n^*$ is the dual (Banach) space of $W_n$.

We also recall that  $T(W_n)$ is  the Banach space

$$\{ \underset{m=0}{\overset{\infty}{\oplus}} u_m: u_m \in \otimes ^m W_n, \sum _{m=0}^{\infty} \|u_m\|<\infty\}$$
endowed with the norm

$$\| \underset{m=0}{\overset{\infty}{\oplus}}  u_m\| = \sum _{m=0}^{\infty} \|u_m\|,$$
and $T_{ \mathbb{C}}(W_n)$ is  the complexification of $T(W_n)$.

Let $H_{F_n}=l^2 (F_n)\otimes H$ and $H_{F_n}\otimes T_{ \mathbb{C}}(W_n)$
 be as
in the definitions of $C^*(F_n, V_n)=C^*(F_n)\otimes C_0(W_n, Cl(W_n))$
and $C^*(F_n, V_n^+)=C^*(F_n)\otimes C_0(W_n, Cl(W_n))^+$.

Given   a collection of elements  $\{w_x\}_{x\in F_n}$ in $T_{ \mathbb{C}}(W_n)$, we define
an operator $A$ on $H_{F_n}\otimes T_{\mathbb{C}} (W_n)$ by:

$$ A((\delta_x \otimes \eta)\otimes \xi) = (\delta_x \otimes \eta)\otimes (w_x \otimes \xi),$$
where $x\in F_n$, $\eta\in H$, $\xi\in T_{ \mathbb{C}}(W_n)$, and $\delta_x$ is the Dirac function at
$x$.

Claim 1: $\|A\| \leq \sup_{x\in F_n} \|w_x\|.$

Proof of  Claim 1:  Let $m\geq 1$.  Given $v\in H_{F_n}\otimes (\otimes^m(W_n))$, we can write
$$v=\sum_{x\in F_{n}} (\delta_x\otimes v_x),$$
where $\delta_x\otimes v_x$ is an element in the closed subspace of $H_{F_n}\otimes (\otimes^m(W_n)) $
spanned by all vectors $(\delta_x\otimes \eta)\otimes \xi$ for all $\eta\in H$ and $\xi\in \,\,\otimes^m(W_n)$.

By the definition of the Banach space norm on $H_{F_n}\otimes T_{\mathbb{C}} (W_n)$, we have

$$\|v\|= \| \sum_{x\in F_n} (\delta_x\otimes v_x)\|
$$ $$=\sup_{ \mu_x\in H^*,\,\, \|\mu_x \|\leq 1,\,\, \lambda_k \in W_n^*,
\,
\|\lambda_k\|\leq 1, \,\,1\leq k\leq m    }\,\,\, ( \sum_{x\in F_n}
  \,\,|(\mu_x\otimes \lambda_1\otimes \cdots \otimes  \lambda_m) (v_x)  |^2)^{\frac{1}{2}}   ,$$
where  $H^*$ is the dual (Hilbert) space of $H$, and $W_n^*$ is the dual (Banach) space of $W_n$.

Claim 1 follows from the above norm formula.

Let $H_{F_n}\otimes Cl(W_n)$ be the Banach space as
in the definitions of $C^*(F_n, V_n)=C^*(F_n)\otimes C_0(W_n, Cl(W_n))$
and $C^*(F_n, V_n^+)=C^*(F_n)\otimes C_0(W_n, Cl(W_n))^+$.

Given   a collection of elements  $\{s_x\}_{x\in F_n}$ in $Cl(W_n)$, we define
an operator $B$ on $H_{F_n}\otimes Cl (W_n)$ by:

$$ B((\delta_x \otimes \eta)\otimes \xi) = (\delta_x \otimes \eta)\otimes (s_x\,\, \xi),$$
where $x\in F_n$, $\eta\in H$, $\xi\in  Cl(W_n)$, and $\delta_x$ is the Dirac function at
$x$, and $s_x\,\, \xi$ is the multiplication of $s_x$ with $\xi$ in $Cl(W_n)$.

Claim 1, together with the definition of the norm on $H_{F_n}\otimes Cl (W_n)$, implies the following:

Claim 2: $\|B\| \leq \sup_{x\in F_n} \|s_x\|.$

For each $w=v\oplus h \in W_n$, let $M_{\phi, w}$ be the bounded linear operator  on $H_{F_n}\otimes Cl(W_n)$
defined by:
$$M_{\phi, w} ((\delta_x \otimes \eta)\otimes \xi) = (\delta_x\otimes \eta ) \otimes (\phi_{x,w} \xi),$$
where $x\in F_n$, $\eta \in H$, $\xi\in Cl(W_n)$, $\delta_x$ is the Dirac function at $x$, $\phi_{x,w}$ is the element in $Cl(W_n)$
defined by: $\phi_{x,w}  = \phi ( (v+f(x))\oplus h)$, and $\phi_{x,w}\xi$ is the multiplication
of $\phi_{x,w}$ with $\xi$ in  $Cl (W_n)$.

By  Claim 2 and the definition of the Banach space norm on $H_{F_n}\otimes Cl(W_n)$,
we have $$\| M_{\phi, w}\| \leq \|\phi\|.$$

Let $C_0 (W_n, H_{F_n}\otimes Cl (W_n))$ be the Banach space in the definitions of
$C^* (F_n, V_n)= C^*(F_n)\otimes C_0(W_n, Cl(W_n))$ and $C^* (F_n, V_n^+)= C^*(F_n)\otimes C_0(W_n, Cl(W_n))^+$ .

We now define an operator $T_{\phi}$ on the Banach space $C_0 (W_n, H_{F_n}\otimes Cl (W_n))$ by:

$$(T_{\phi} \zeta)(w)= M_{\phi, w} (\zeta(w))$$
for all $\zeta \in C_0 (W_n, H_{F_n}\otimes Cl (W_n))$ and $w\in W_n$.

By the above inequality, we have $$\| T_{\phi}\| \leq \|\phi\|.$$

Let $\psi(b_n)$ be the operator on the Banach space $C_0 (W_n, H_{F_n}\otimes Cl (W_n))$ defined by:
$$ ( \psi(b_n) \zeta)(w)= (b_n\otimes 1) \zeta(w)$$
for all $\zeta \in C_0 (W_n, H_{F_n}\otimes Cl (W_n))$ and $w\in W_n$.

By the definition of the norm on the Banach space $C_0 (W_n, H_{F_n}\otimes Cl (W_n))$,
we have

$$ \| \psi(b_n)\|\leq \|b_n\|.$$

Note that $c_n (b_n,\phi)$ is the product of $T_{\phi}$ with $\psi(b_n)$ as operators on
the Banach space $C_0 (W_n,  H_{F_n}\otimes  Cl(W_n))$.

It follows that $$\| c_n(b_n, \phi)\|\leq \|b_n\|  \,\, \|\phi\|$$
for all $\phi\in C_0(W_n, Cl(W_n))^+$.

This proves  part (1) of Lemma 6.1.

Recall that $index(F_{V_n, \varphi})$ is an element in $C_0(W_n, Cl(W_n))^+ \otimes M_2 (\mathbb{C})$.
By part (2) of Proposition 3.2, there exists $C>0$ such that $index(F_{V_n, \varphi})\leq C$ for all $n$.
Now part (2) of  Lemma 6.1 follows from part (1) of Lemma 6.1,  Proposition 3.3, and the bounded geometry property of $\Gamma$.
\begin{flushright}
$\Box$
\end{flushright}

Define \ $ q_0(e')=\underset{n=1}{\overset{\infty}{\oplus}} b_n
\in C^*_{alg}(\left\{F_n,V^+_n\right\})\otimes M_2 (\mathbb{C})$
by:
\[
(b_n(x,y))(v\oplus h)=e'_n (x,y)\otimes \left(\begin{array}{cc}1&0
\\ 0&0\end{array}\right)\,
\]
for all $x,y\in F_n$, and $ v\oplus h \in W_n=V_n\oplus V^*_n$.

Let $\varepsilon$  be as in Proposition 6.1. Choose $\varepsilon$
to be sufficiently small.  Let $p(e')$ and $q(e')$ be idempotents
in $C^* (\{F_n, V_n^+\})\otimes M_2 (\mathbb{C})$ obtained from
$p_0(e')$ and $q_0(e')$ by functional calculus. Note that
$p(e')-q(e')\in C^*(\left\{F_n,V_n\right\})\otimes M_2
(\mathbb{C})$.

\begin{defi}
We define $\beta_{\{V_n\}}[e]\in K_0(C^*(\left\{F_n,V_n\right\}))$
to be the difference element $D(p(e'),q(e'))$.
\end{defi}

By suspension, we can similarly define
\[
\beta_{\left\{{V_n}\right\}}\,:\,K_1(C^*(\left\{F_n\right\})) \to K_1
(C^*(\left\{F_n,V_n\right\})).
\]

Notice that the homomorphism
$\beta_{\left\{{V_n}\right\}}\,:\,K_*(C^*(\left\{F_n\right\})) \to
K_* (C^*(\left\{F_n,V_n\right\})),$ is independent of the choice
of $\varphi.$

We can  define the following Bott map by the same method:

\[
\beta_{\left\{{V_n}\right\}}:K_*(C^*(\left\{F_n,\partial_{\Gamma}
F_n\right\}))\to K_* (C^*(\left\{F_n,\partial_{\Gamma}
F_n,V_n\right\}_n)).
\]

For any $R>0$, we define $C^*(F_n,V_n)_R$ to be the closed
subalgebra of $C^*(F_n,V_n)$ generated by elements $a\in
C^*(F_n,V_n)$ such that $propagation\,(a)<R$ and
$support\,(a(x,y))\subseteq B_{W_n}(f(x)\oplus 0,R)$ for all
$x,y\in F_n$. Note that there exists $R'>0$ such that every
element in $C^*(F_n,V_n)_R$ has propagation at most $R'$, where $R'$ does not depend on $n$.
\begin{defi}
\[
C^*(\left\{F_n,V_n\right\}_n)_R=\left\{\underset{n=1}{\overset{\infty}{\oplus}}
a_n\,:\,a_n \in C^*(F_n,V_n)_R\,,\
\sup_n\|a_n\|<+\infty\right\}\,,
\]
where $C^*(\left\{F_n,V_n\right\}_n)_R$ is endowed with the norm:
$\left\|\underset{n=1}{\overset{\infty}{\oplus}} a_n\right\|=
\sup_n\|a_n\|$.
\end{defi}

We can similarly define the Banach algebras
$\,\,\,\,\,C^*(\left\{\partial_\Gamma
F_n,V_n\right\}_n)_R\,\,\,\,\,$ and
$\,\,\,\,\,C^*(\{F_n,\pt_\Gamma F_n,V_n\}_n)_R\,.$
\begin{lem} We have
\begin{enumerate}
\item[(1)] $C^*(\left\{F_n,V_n\right\}_n)=\lim_{R\to\infty}
C^*(\left\{F_n,V_n\right\}_n)_R$\,;
\item[(2)]$C^*(\left\{\partial_\Gamma F_n,V_n\right\}_n)=\lim_{R\to\infty}
C^*(\left\{\partial_\Gamma F_n,V_n\right\}_n)_R$\,;
\item[(3)]$C^*(\left\{F_n,\partial_\Gamma F_n,V_n\right\}_n)=\lim_{R\to\infty} C^*(\left\{F_n
,\partial_\Gamma F_n,V_n\right\}_n)_R$\,.
\end{enumerate}
\end{lem}

Notice that there exists a natural homomorphism:
$$\lim_{R\to\infty} K_*(C^*(\left\{F_n,V_n\right\}_n)_R)
\to\lim_{R\to\infty}\underset{n=1}{\overset{\infty}{\oplus}}K_*
(C^*(F_n,V_n)_R)\,,$$
where $\underset{n=1}{\overset{\infty}{\oplus}}K_*(C^*(F_n,V_n)_R)=\{
\underset{n=1}{\overset{\infty}{\oplus}}z_n:\,\, z_n\in K_*(C^*(F_n,V_n)_R)\}.$

Similarly there exists a natural homomorphism:
$$\lim_{R\to\infty} K_*(C^*(\left\{F_n ,\partial_\Gamma
F_n,V_n\right\}_n)_R)\to
\lim_{R\to\infty}\underset{n=1}{\overset{\infty}{\oplus}} K_*(
C^*(F_n,\partial_\Gamma F_n,V_n)_R)\,,$$
where $\underset{n=1}{\overset{\infty}{\oplus}}K_*(C^*(F_n,\partial_\Gamma F_n, V_n)_R)=\{
\underset{n=1}{\overset{\infty}{\oplus}}z_n:\,\, z_n\in K_*(C^*(F_n,\partial_\Gamma F_n, V_n)_R)\}.$

Composing $ \beta_{\left\{V_n\right\}}$ with the above
homomorphisms, we obtain homomorphisms (still denoted by
$\beta_{\left\{V_n\right\}})$:
\[
K_*(C^*(\left\{F_n\right\}_n))\to\lim_{R\to\infty}\underset{n=1}{\overset{\infty}{\oplus}}K_*
(C^*(F_n,V_n)_R)\,,
\]

\[
K_*(C^*(\left\{F_n,\partial_\Gamma F_n\right\}_n))\to
\lim_{R\to\infty}\underset{n=1}{\overset{\infty}{\oplus}} K_*(
C^*(F_n,\partial_\Gamma F_n,V_n)_R)\,.
\]

\begin{prop} Assume that $\{V'_n\}$ is a sequence of finite
dimensional subspaces of $X$ such that $V_n \subseteq V'_n$. For any $[z]\in K_*
( C^*(\left\{F_n,\partial_\Gamma F_n\right\}_n))$,\\
$$\beta_{\left\{{V_n}\right\}}[z]=0\,\,\, \mathrm{ in}\,\,\,
\lim_{R\to\infty}\underset{n=1}{\overset{\infty}{\oplus}} K_*(
C^*(F_n,\partial_\Gamma F_n,V_n)_R)$$ if and only if
$$\beta_{\left\{{V_n'}\right\}}[z]=0\,\,\,
\mathrm{in} \,\,\,
\lim_{R\to\infty}\underset{n=1}{\overset{\infty}{\oplus}} K_*(
C^*(F_n,\partial_\Gamma F_n,V_n')_R).$$
\end{prop}

Roughly speaking, the proof of Proposition 6.5 is based on the
idea that the product of  Bott elements for two finite dimensional
vector spaces is  the Bott element for the product space of the
two vector spaces. We need some preparations before we can  prove
Proposition 6.5.

Next we shall introduce the concept of a Fredholm pair for Banach algebras
and define its index.
The concept of Fredholm pair  and its index is motivated by KK-theory. In the special case of $C^*$-algebras,
it is compatible with the corresponding construction in KK-theory.

Let $B$ be a graded  unital complex Banach algebra with the
grading induced by a grading operator $\varepsilon$ in $B$
satisfying $\varepsilon^2 =1$ and $\|\varepsilon\|=1$. Let $A$ be
a graded closed two sided ideal in $B$.

Let $F$ be an element of degree 1 in $B $ and let $e$ be an
idempotent of degree 0 in $B$.

$(F, e)$ is said to be a Fredholm pair if
\begin{enumerate}
\item[(1)] $e(F^2-1)\in A\,;$
\item[(2)] $eF-Fe \in A$\,.
\end{enumerate}

We define $Index(F,e)\in K_0(A)$  of a Fredholm pair $(F,e)$ as
follows.

Let
$$ p=\left(\frac{1+\varepsilon}{2}\right) e
\left(\frac{1+\varepsilon}{2}\right),$$
$$ q=\left(\frac{1-\varepsilon}{2}\right) e
\left(\frac{1-\varepsilon}{2}\right),$$
$$ \alpha=\left(\frac{1+\varepsilon}{2}\right) F
\left(\frac{1-\varepsilon}{2}\right),$$
$$ \alpha'=\left(\frac{1-\varepsilon}{2}\right) F
\left(\frac{1+\varepsilon}{2}\right).$$

Define
\[
a_{11}(F,e)=1+(p-p \alpha q\alpha'p)p\alpha p\alpha'p +(p\alpha
q\alpha'p-p),
\]
\[
a_{12}(F,e)=(p-p \alpha q\alpha'p)p\alpha q(q-q \alpha'p \alpha
q)+ p\alpha q (q-q\alpha p \alpha'q),
\]
\[
a_{21}(F,e)=(q-q \alpha'p \alpha q)q\alpha' p,
\]
\[
a_{22}(F,e)=(q-q \alpha'p \alpha q)^2.
\]
Define $$index(F,e)=\left(\begin{array}{cc}a_{11}(F,e)&a_{12}(F,e)\\
a_{21}(F,e)& a_{22}(F,e)
\end{array}\right).$$

It is not difficult to verify that $index(F,e)$ is an idempotent
in $M_2(A)$.

We define
\[
Index(F,e)=[index(F,e)]-\left[\left(\begin{array}{cc}1&0\\0&0
\end{array}\right)\right]\in K_0(A),
\]
where $K_0(A)$ is the K-group of $A$ considered as a Banach
algebra without grading.

We remark that if $e=1,$ then
$$Index (F,e)=Index (F),$$ where $Index (F)$ is defined as in
Section 2.

Let $A_1$ and $A_2$ be two Banach algebras. Assume that  $A_1\otimes A_2$ is  a Banach algebra tensor product of $A_1$
and $A_2$.
We shall need an
explicit construction of the product:
\[
K_0(A_1)\times K_0(A_2)\to K_0(A_1\otimes A_2).
\]

Let $\pi_i:A^+_i\to\mathbb{C}$ be the homomorphism defined by
$\pi_i(a+cI)=c$ for any $a\in A_i,\,c\in \mathbb{C}$ and $i=1,2$,
where $A_i^+$ is obtained from $A_i$ by adjoining the identity.

Let $\pi:A^+_1\otimes A_2 + A_1 \otimes A^+_2 \to \mathbb{C}
\otimes A_2 + A_1 \otimes \mathbb{C}$ be the homomorphism defined by
\[
\pi (a_1\otimes a_2 + b_1\otimes b_2)=\pi_1 (a_1)\otimes a_2 +
b_1\otimes \pi_2 ( b_2)
\]
for any $a_1\otimes a_2\in A^+_1 \otimes A_2,\ b_1\otimes b_2 \in
A_1 \otimes A^+_2 $. We have the following split exact sequence:
\[
0 \to A_1 \otimes A_2 \to A^+_1 \otimes A_2 +A_1 \otimes
 A^+_2 \underset{i}{\overset{\pi}{\rightleftarrows}} \mathbb{C}
\otimes A_2 + A_1 \otimes \mathbb{C}\to 0,
\]
where $i:\mathbb{C} \otimes A_2 + A_1 \otimes \mathbb{C}\to A^+_1
\otimes A_2 +A_1 \otimes A^+_2$ is the inclusion homomorphism.

Next we shall construct an explicit homomorphism
$$\phi\,:\,Ker\, \pi_*
\to K_0(A_1\otimes A_2).$$

Let $p$ be an idempotent in $M_n((A_1^+ \otimes A_2 + A_1\otimes
A_2^+)^+)$ such that $[p]-[i\pi(p)]$ represents an element in
$Ker\, \pi_*$.

Let
\[
Z=\left(\begin{array}{cccc}i\pi(p)&0&1-i\pi(p)&0 \\ 1-i
\pi(p)&0&0&i\pi(p)\\0&0&i\pi(p)&1-i\pi(p)\\0&1&0&0\end{array}\right)\,.
\]

We have
\[
Z^{-1}=\left(\begin{array}{cccc}i\pi(p)&1-i\pi(p)&0&0
\\0&0&0&1\\1-i\pi(p)&0& i\pi(p)&0\\0&i\pi(p)&1-i\pi(p)&0\end{array}\right)\,.
\]

Let
\[
\phi_0(p)=Z^{-1}\left(\begin{array}{cccc}p&0&0&0
\\0&1-i\pi(p)&0&0\\0&0&0&0\\0&0&0&0\end{array}\right)Z\,.
\]

We define
\[
\phi\ ([p]-[i\pi(p)])=[\phi_0(p)]-\left[\left(\begin{array}{cccc}1&0&0&0\\0&0&0&0\\0&0&0&0\\0&0&0&0
\end{array}\right)\right]\in K_0(A_1\otimes A_2).
\]

\begin{lem} The homomorphism $\phi$ from $Ker\, \pi_*$ to $ K_0(A_1\otimes A_2) $ is an isomorphism.
\end{lem}

The product: $ K_0 (A_1)\times  K_0 (A_2)\to  K_0 (A_1 \otimes A_2)$
can now be described as follows. Given idempotents $p_0,p_1$\ in $M_n(A^+_1)$ and $q_0,q_1$ in $M_n
(A^+_2)$ representing $[p_0]-[p_1]\in K_0 (A_1)$ and $[q_0]-[q_1]
\in K_0 (A_2)$ such that $\pi_1(p_0)=\pi_1(p_1),\ \pi_2(q_0)=\pi_2(q_1)$, note that
\[
D((p_0\otimes q_0)\oplus(p_1\otimes q_1),\,(p_0\otimes
q_1)\oplus(p_1\otimes q_0))\in Ker\,\pi_*\subseteq
K_0(A_1^+\otimes A_2+A_1\otimes A_2^+)\,.
\]
We define the product of $[p_0]-[p_1]$ and $[q_0]-[q_1]$ to be
\[
\phi(D((p_0\otimes q_0)\oplus(p_1\otimes q_1),\,(p_0\otimes
q_1)\oplus(p_1\otimes q_0)))\in K_0(A_1\otimes A_2)\,.
\]

Let $B_1$ and $B_2$ be two  graded unital complex  Banach algebras
with gradings induced by grading operators $\varepsilon_1$ and
$\varepsilon_2$ (respectively in $B_1$ and $B_2$) satisfying
$\varepsilon_i^2=1$ and $\| \varepsilon_i\|=1$ for $i=1,2$. Let
$A_1$ and $A_2$ be respectively graded closed two sided ideals in
$B_1$ and $B_2$.

Assume that $F_i$ is an element of degree 1 in $B_i$ satisfying
\begin{enumerate}
\item[(1)] $F_i^2-1\in A_i$\,,
\item[(2)] $(1-F_i^2)^{1/2}$ is a well defined element in
$A_i$,
\end{enumerate}
where $i=1,\,2$.

Assume that $B_1 \otimes B_2$ is a Banach tensor product of
$B_1$ and $B_2$.
Endow $B_1\otimes B_2$ with the grading induced by the grading
operator $\varepsilon_1 \otimes\varepsilon_2$.

Let $M$ and $N$ be elements in $B_1 \otimes B_2$ such that
\begin{enumerate}
\item[(1)] $M$ and $N$ have degree $0$ and respectively commute with $F_1
\otimes 1$ and $1\otimes F_2$;
\item[(2)] $M^2 +N^2-1 \in A_1 \otimes A_2;$
\item[(3)] $M((1-F_1^2)^{1/2}\otimes 1)$ and $N(1\otimes (1-F_2^2)^{1/2}) $ are elements in $
A_1\otimes A_2.$
\end{enumerate}

Define
\[
F=M(F_1\otimes 1)+ N(\varepsilon_1 \otimes F_2)\in B_1\otimes B_2.
\]
It is easy to check that $F$ has degree 1 and $ F^2-1\in
A_1\otimes A_2.$

Hence we can define
\[
Index(F)\in K_0(A_1\otimes A_2).
\]

\begin{prop}
Assume that there exist homotopies $M_t$ and $N_t$ ($t\in [0,1]$)
in $B_1\otimes B_2$ such that
\begin{enumerate}
\item[(1)] $M_t$ and $N_t$ have degree 0 for each $t\in [0,1]$;
\item[(2)] $M_t$ and $N_t$ commute respectively with $F_1\otimes 1$ and $1\otimes
F_2$ for each $t\in [0,1]$;
\item[(3)] $M_t$ and $N_t$ commute respectively with $(1-F_1^2)^{1/2}\otimes 1$ and $1\otimes
(1-F_2^2)^{1/2}$ for each $t\in [0,1]$;
\item[(4)] $M_t^2+N_t^2-1 \in A_1\otimes A_2$ for all
$t\in [0,1]$;
\item[(5)] $M_t ((A_1^+\otimes A_2)\oplus (A_1\otimes A_2^+)
)\subseteq$ $(A_1^+\otimes A_2)\oplus (A_1\otimes A_2^+) $ for all
$t\in [0,1]$,

{\noindent $((A_1^+\otimes A_2)\oplus (A_1\otimes A_2^+) )M_t
\subseteq$ $(A_1^+\otimes A_2)\oplus (A_1\otimes A_2^+) $ for all
$t\in [0,1]$,}

 {\noindent $N_t
((A_1^+\otimes A_2)\oplus (A_1\otimes A_2^+) )\subseteq $
$(A_1^+\otimes A_2)\oplus (A_1\otimes A_2^+) $ for all $t\in
[0,1]$,}

{\noindent and $((A_1^+\otimes A_2)\oplus (A_1\otimes A_2^+) ) N_t
\subseteq $ $(A_1^+\otimes A_2)\oplus (A_1\otimes A_2^+) $ for all
$t\in [0,1]$;}
\item[(6)] $M_0=M,\, N_0=N,\,$ $M_1=N_1=1/\sqrt{2}$.
\end{enumerate}
Then $Index(F)$ is equal to the product of $Index(F_1)$ with
$Index(F_2)$ in

{\noindent $K_0(A_1\otimes A_2)$.}
\end{prop}
{\noindent {\it Proof:}}

Let \[ \tilde{F_i}=\left(\begin{array}{cc}F_i& 0\\
0& 0
\end{array}\right)\,, \quad i=1,\,2\, ,
\]

\[ p_0=\left(\begin{array}{cc}1 & 0\\
0& 0
\end{array}\right)\,.
\]

For $i=1,2,$ let $M_2 (B_i)$ be endowed with the grading induced
by the grading operator:
\[ \left(\begin{array}{cc}\varepsilon_i & 0\\
0& -\varepsilon_i
\end{array}\right)\,.
\]

$(\tilde{F_i}, p_0)$ is a Fredholm pair for $i=1,2$. We have
$$Index(\tilde{F_i}, p_0)=Index (F_i).$$

Let
\[ W_i=\frac{1}{2} \left(\begin{array}{cc}1+\varepsilon_i& 1-\varepsilon_i\\
1-\varepsilon_i& 1+\varepsilon_i
\end{array}\right)\,, \quad i=1,\,2\,,
\]

\[\varepsilon_0=
\left(\begin{array}{cc}1 & 0\\
0& -1
\end{array}\right).\]

We have
\[ W_i^{-1} \left(\begin{array}{cc}\varepsilon_i & 0\\
0& -\varepsilon_i
\end{array}\right)W_i=   \varepsilon_0 \,, \quad i=1,\,2\,.
\]

Let $$E_i=W_i^{-1}\tilde{F_i} W_i, \,\,\,e_i=W_i^{-1} p_0 W_i\,$$
for $i=1,\,2\,.$

Notice that $(E_i, e_i)$ is a Fredholm pair with respect to the
grading on $M_2 (B_i)$ induced by the grading operator
$\varepsilon_0$ for $i=1,2$. We have
$$Index (E_i, e_i)= Index (F_i)$$ for
$i=1,\,2\,.$

Let
\[
G_i=\left(\begin{array}{cc}E_i& (1-E_i^2)^{1/2}\\ (1-E_i^2)^{1/2}&
-E_i
\end{array}\right)\,, \quad i=1,\,2\,.
\]

Define
\[
F'=M'(G_1\otimes 1)+ N'((\varepsilon_0 \oplus
-\varepsilon_0)\otimes G_2)) \in M_2(M_2(B_1))\otimes
M_2(M_2(B_2))\,,
\]
where $\,\,\,\,\,\,M'=\sum_{k=1}^{16} M,\,\,\,\,\,\,\,\,
N'=\sum_{k=1}^{16} N\,\,\,\,\,\,\in \,\,\,\,M_{16} (B_1\otimes
B_2) \,\,\,\,\cong \,\,\,\,\,\,\,\,\,\,\,\,\,\,\,\,\,$\break $
M_2(M_2(B_1))\otimes M_2 (M_2(B_2)).$

Endow $M_2(M_2(B_1))\otimes M_2(M_2(B_2))$ with the grading
operator $$\varepsilon'=(\varepsilon_0\oplus
-\varepsilon_0)\otimes (\varepsilon_0\oplus -\varepsilon_0).$$

 Define
\[
e'=\left(\begin{array}{cc}e_1& 0\\ 0 & 0
\end{array}\right)\otimes \left(\begin{array}{cc}e_2& 0\\ 0 & 0
\end{array}\right)\in M_2(M_2(B_1))\otimes M_2(M_2(B_2))\,.
\]

It is not difficult to verify that $(F', e')$ is a Fredholm pair
with respect to the grading given by the grading operator
$\varepsilon'$. We have
$$Index(F)=Index(F',e').$$

Define
\begin{eqnarray*}
F''&=&M''((G_1\oplus -G_1)\otimes 1)+N''(( (\varepsilon_0\oplus
-\varepsilon_0)\oplus(-\varepsilon_0\oplus
\varepsilon_0))\otimes(G_2\oplus -G_2))
\\ &\in& M_2(M_4(B_1))\otimes M_2(M_4(B_2))\,,
\end{eqnarray*}
where $M''= M'\oplus M' \oplus M'\oplus M',\, N''=N' \oplus
N'\oplus N' \oplus N'$, and

{\noindent $M_2(M_4(B_1))\otimes M_2(M_4(B_2))$ is endowed with
the grading operator
\[
\varepsilon''=((\varepsilon_0\oplus
-\varepsilon_0)\oplus(-\varepsilon_0\oplus \varepsilon_0))\otimes
((\varepsilon_0\oplus -\varepsilon_0)\oplus(-\varepsilon_0\oplus
\varepsilon_0))\,.
\] }

Define
\begin{eqnarray*}
e''&=& \left(\left(\begin{array}{cc}e_1& 0\\ 0 & 0
\end{array}\right)\oplus \left(\begin{array}{cc}0 & 0\\ 0 & 0
\end{array}\right)\right)\otimes \left(\left(\begin{array}{cc}e_2& 0\\ 0 & 0
\end{array}\right)\oplus \left(\begin{array}{cc}0 & 0\\ 0 & 0
\end{array}\right)\right)
\\ & & \in M_2(M_4(B_1))\otimes M_2(M_4(B_2))\,.
\end{eqnarray*}

Notice that $(F'', e'')$ is a Fredholm pair with respect to the
grading given by the grading operator $\varepsilon''$.  We have
$$Index(F',e')=Index(F'',e'').$$

Let
\[
U_i=\frac{1}{\sqrt{2}}\left(\begin{array}{cc}1& -G_i\\ G_i & 1
\end{array}\right)
\]
for $i=1,2$.

 We have
\[
U_i\left(\begin{array}{cc}G_i& 0\\ 0 & -G_i
\end{array}\right)U_i^{-1}=\left(\begin{array}{cc}0& 1\\ 1 & 0
\end{array}\right)\,
\]
for $i=1,2$.

 Notice that $U_1\otimes 1$ and $1\otimes U_2$ have
degree 0 in $M_2(M_4(B_1))\otimes M_2(M_4(B_2))$. Let
\begin{eqnarray*}
F'''&=& (1\otimes U_2)(U_1\otimes 1)F''(U_1\otimes
1)^{-1}(1\otimes U_2)^{-1}
\\ &=&M''\left(\left(\begin{array}{cc}0& 1\\ 1 & 0
\end{array}\right)\otimes 1\right) + N''\left(((\varepsilon_0\oplus -\varepsilon_0)\oplus(-\varepsilon_0\oplus
\varepsilon_0))\otimes \left(\begin{array}{cc}0& 1\\ 1 & 0
\end{array}\right)\right)\,,
\end{eqnarray*}

\begin{eqnarray*}
e'''= \left(U_1\left(\left(\begin{array}{cc}e_1& 0\\ 0 & 0
\end{array}\right)\oplus \left(\begin{array}{cc}0& 0\\ 0 & 0
\end{array}\right)\right)U_1^{-1}\right)\otimes \left(U_2\left(\left(\begin{array}{cc}e_2& 0\\ 0 & 0
\end{array}\right)\oplus \left(\begin{array}{cc}0& 0\\ 0 & 0
\end{array}\right)\right)U_2^{-1}\right)\,.
\end{eqnarray*}

Notice that $\varepsilon''$ commutes with  $(1\otimes
U_2)(U_1\otimes 1)$. It follows that $(F''', e''')$ is a Fredholm
pair with respect to the grading induced by the grading operator
$\varepsilon''$. We have
$$Index(F''',e''')=Index(F'',e'').$$

Let $$F_i'''=\left(\begin{array}{cc}0& 1\\ 1 & 0
\end{array}\right) \in M_2 (M_4(B_i)),$$
\[
e_i'''=U_i\left(\left(\begin{array}{cc}1& 0\\ 0 & 0
\end{array}\right)\oplus \left(\begin{array}{cc}0& 0\\ 0 & 0
\end{array}\right)\right)U_i^{-1} \in M_2 (M_4(B_i)) \,,\quad i=1,\,2.
\]

Notice that $(F_i''', e_i''')$ is a Fredholm pair with respect to
the grading given by the grading operator $ (\varepsilon_0\oplus
-\varepsilon_0)\oplus(-\varepsilon_0\oplus \varepsilon_0)$ in $M_2
(M_4(B_i))$.

We have
$$Index(F_i)=Index(F_i''',e_i''')$$ for $i=1,2$.

Let $$ V=\frac{1}{2}\left(\begin{array}{cc}1+(\varepsilon_0\oplus -\varepsilon_0)& 1-(\varepsilon_0\oplus -\varepsilon_0)\\
1-(\varepsilon_0 \oplus -\varepsilon_0)&
1+(\varepsilon_0\oplus-\varepsilon_0)
\end{array}\right).$$
We have
$$ V^{-1}\left(\begin{array}{cc}\varepsilon_0\oplus -\varepsilon_0& 0\\
0 & -\varepsilon_0\oplus \varepsilon_0
\end{array}\right)V=   \left(\begin{array}{cc}1& 0\\
0 & -1
\end{array}\right) .$$
Note also that $V$ commute with $\left(\begin{array}{cc}0& 1\\
1 & 0
\end{array}\right).$

Write $$V^{-1}e_i '''V=\left(\begin{array}{cc}p_i& 0\\ 0 & q_i
\end{array}\right)\in M_2 (M_4(B_i)) \,,\quad i=1,\,2.$$

Let
$$e_i''''=\left(\begin{array}{cccc}
  p_i & 0 & 0 & 0 \\
  0 & 1-q_i & 0 & 0 \\
  0 & 0 & 0& 0 \\
  0 & 0 & 0 & 0
\end{array}\right) \oplus
\left(\begin{array}{cccc}
  q_i & 0 & 0 & 0 \\
  0 & 1-q_i & 0 & 0 \\
  0 & 0 & 0& 0 \\
  0 & 0 & 0 & 0
\end{array}\right) \in M_{2\times 4^2}(B_i),$$
$$F''''=M''' \left( \left(\begin{array}{cc}
  0 & I' \\
  I' & 0
\end{array}\right)\otimes I''\right) + N'''\left(\tau\otimes
 \left(\begin{array}{cc}
  0 & I' \\
  I' & 0
\end{array}\right)\right) \in M_{2\times 4^2}(B_1)\otimes M_{2\times 4^2}(B_2),$$
where $I'=I\oplus I\oplus I \oplus I$ ($I$ is the identity element
in $M_4(B_1)$ or $M_4(B_2)$), $I''=I'\oplus I'$, $M'''=M''\oplus
M''\oplus M''\oplus M''$, $N'''=N''\oplus N''\oplus N''\oplus
N''$, and
$$\tau= \left(\begin{array}{cc}
  I' & 0 \\
  0 & -I'
\end{array}\right).$$

Let $$e''''=e_1''''\otimes e_2''''\in M_{2\times 4^2}(B_1)\otimes
M_{2\times 4^2}(B_2).$$

Observe that $(F'''', e'''')$ is a Fredholm pair with respect to
the grading induced by the grading operator $\tau\otimes \tau$. We
have
$$Index(F''',e''')= Index (F'''', e'''').$$

Let $$Z(p_i, q_i)=\left(\begin{array}{cccc}
  q_i & 0 & 1-q_i & 0 \\
  1-q_i & 0 & 0 & q_i \\
  0 & 0 & q_i & 1-q_i \\
  0 & 1 & 0 & 0
\end{array}\right).$$
Write
$$ (Z(p_i,q_i)\oplus Z(p_i, q_i))^{-1} e_i'''' (Z(p_i,q_i)\oplus Z(p_i,
q_i))= \left( \begin{array}{cc}
  p_i' & 0 \\
  0 & q_i'
\end{array}\right) \in M_2(M_{4^2}(A_i)).$$

We have $$p_i'-q_i'\in M_{4^2}(A_i)\,\,\,\,\, \, \,\,\,\, \cdots
\cdots \cdots \,\,\,\,\,\,\,\,\,\,(J).$$

Notice that $(Z(p_1,q_1)\otimes 1) \oplus (Z(p_1, q_1)\otimes 1)$
and $(1\otimes Z(p_2,q_2))\oplus (1\otimes Z(p_2,q_2))$ commute
with $F''''.$

Define
$$ H(t)= M_t' \left( \left(\begin{array}{cc}
  0 & I' \\
  I' & 0
\end{array}\right)\otimes I'' \right) + N_t'\left(\tau\otimes
 \left(\begin{array}{cc}
  0 & I' \\
  I' & 0
\end{array}\right)\right) \in M_{2\times 4^2}(B_1)\otimes M_{2\times 4^2}(B_2),$$
where $M_t'=\underset{k=1}{\overset{2\times 4^4}{\oplus}} M_t$ and
$N_t'=\underset{k=1}{\overset{2\times 4^4}{\oplus}} N_t$.

Let $$g= \left( \begin{array}{cc}
  p_1' & 0 \\
  0 & q_1'
\end{array}\right) \otimes
\left( \begin{array}{cc}
  p_2' & 0 \\
  0 & q_2'
\end{array}\right) \in M_{2\times 4^2}(B_1)\otimes M_{2\times 4^2}(B_2).$$

We have
$$H(t)g -g H(t)\in (M_{2\times 4^2}(A_1^+)\otimes M_{2\times 4^2}(A_2)) +
(M_{2\times 4^2}(A_1)\otimes M_{2\times 4^2}(A_2^+)),$$
$$ H(t)^2 -1 \in M_{2\times 4^2}(A_1)\otimes M_{2\times 4^2}(A_2).$$

Hence, by the formula of $index (H(t),g)$, condition (5) of the
proposition and the above property $(J)$ of $p_i'$ and $q_i'$, we
know that $index (H(t), g)$ is a homotopy  of idempotents in
$M_2(((M_{2\times 4^2}(A_1^+)\otimes M_{2\times 4^2}(A_2)) +
(M_{2\times 4^2}(A_1)\otimes M_{2\times 4^2}(A_2^+)))^+)$ (note
that although $(H(t), g)$ may  not be a Fredholm pair, $index
(H(t),g)$ is still well defined as an idempotent for each $t$).

It follows that
$$Index(H(1), g)=
Index(H(0), g)$$
 in $K_0 ((M_{2\times 4^2}(A_1^+)\otimes M_{2\times 4^2}(A_2)) +
(M_{2\times 4^2}(A_1)\otimes M_{2\times 4^2}(A_2^+))).$

This implies that
$$Index
(H(1), g)\in Ker \, \pi_* ,$$
 where $\pi$ is the natural
homomorphism from $(M_{2\times 4^2}(A_1^+)\otimes M_{2\times
4^2}(A_2)) + (M_{2\times 4^2}(A_1)\otimes M_{2\times 4^2}(A_2^+))$
to $(M_{2\times 4^2}(\mathbb{C})\otimes M_{2\times 4^2}(A_2)) +
(M_{2\times 4^2}(A_1)\otimes M_{2\times 4^2}(\mathbb{C})).$

 Notice that
$$ H(1)= \sqrt{\frac{1}{2}} \left(\left( \left(\begin{array}{cc}
  0 & I' \\
  I' & 0
\end{array}\right)\otimes I'' \right) + \left(\tau\otimes
 \left(\begin{array}{cc}
  0 & I' \\
  I' & 0
\end{array}\right)\right)\right).$$

It is not difficult to see that there exists a homotopy $H_1 (t)$
in

{\noindent $M_{2\times 4^2}(\mathbb{C} )\otimes M_{2\times 4^2}(
\mathbb{C} )\subseteq M_{2\times 4^2}(B_1)\otimes M_{2\times
4^2}(B_2)$
 such that}
\begin{enumerate}
\item[(1)] $H_1 (t)$ has degree one for all $t\in [0,1]$;
\item[(2)] $ H_1 (0)= H(1)$;
\item[(3)] $ H_1 (1)= \left(\begin{array}{cc}
  0 & 1 \\
  1 & 0
\end{array}\right) \in$

$M_2 (M_2((M_{2\times 4}(B_1)\otimes M_{2\times 4}(B_2)))) ,$
 where $M_2 (M_2((M_{2\times 4}(B_1)\otimes M_{2\times 4}(B_2)) ))$ is a
graded algebra with the grading operator $\left(\begin{array}{cc}
  1 & 0 \\
  0 & -1
\end{array}\right)$ and is identified with
$M_{2\times 4^2}(B_1)\otimes M_{2\times 4^2}(B_2)$ as graded
algebras;

\item[(4)] $(H_1 (t))^2=\left(\begin{array}{cc}
  1 & 0 \\
  0 & 1
\end{array}\right);$
\item[(5)] $H_1 (t)g -gH_1 (t)\in$

$M_2 (M_2((M_{2\times 4}(A_1^+)\otimes M_{2\times 4}(A_2)) +(
M_{2\times 4}(A_1)\otimes M_{2\times 4}(A_2^+))))\subseteq$

 $M_2 (M_2((M_{2\times 4}(B_1)\otimes M_{2\times 4}(B_2)) ))
$ for all $t\in [0,1].$
\end{enumerate}

By the formula of $index(H_1 (t), g)$ and the above properties of
$H_1(t)$, we know that $index (H_1 (t), g)$ is a homotopy of
idempotents in

{\noindent $M_2(((M_{2\times 4^2}(A_1^+)\otimes M_{2\times
4^2}(A_2)) + (M_{2\times 4^2}(A_1)\otimes M_{2\times
4^2}(A_2^+)))^+)$} (note that although $(H_1(t), g)$ may not be a
Fredholm pair, $index (H_1(t),g)$ is still well defined as an
idempotent for each $t$).

It follows that
$$Index(H_1(1), g)= Index(H_1 (0), g)=j_*(Index (F))$$ in $K_0
((M_{2\times 4^2}(A_1^+)\otimes M_{2\times 4^2}(A_2)) +
(M_{2\times 4^2}(A_1)\otimes M_{2\times 4^2}(A_2^+))),$ where $j$
is the natural inclusion homomorphism from $A_1\otimes A_2$ to
$(M_{2\times 4^2}(A_1^+)\otimes M_{2\times 4^2}(A_2)) +
(M_{2\times 4^2}(A_1)\otimes M_{2\times 4^2}(A_2^+)).$

 Now
Proposition 6.7 follows from the fact that $\phi(Index(H_1(1),g))$
is the product of $Index(F_1''',e_1''')$ with
$Index(F_2''',e_2''')$, where $\phi$ is as in Lemma 6.6.

\begin{flushright}
$\Box$
\end{flushright}

Let $V_n$ and $V_n'$ be as in Proposition 6.5. Let $V_n''$ be a
finite dimensional subspace of $X$ such that $V_n'= V_n\oplus
V_n''$.

 Define
$$\|x\|_1=\sqrt{\|x_1\|^2+\|x_2\|^2}$$
$$x^{*,1}= x_1^* \oplus x_2^*$$
 for any $x=x_1\oplus x_2\in V_n\oplus V_n''$, and
$$\|h\|_1=\sqrt{\|h_1\|^2+\|h_2\|^2}$$
$$h^{*,1} =h_1^* \oplus h_2^* $$
for any $h=h_1\oplus h_2\in (V_n')^*=V_n^*\oplus (V_n'')^*$.

We define $F_{V_n', \varphi,1}
\in C_b (W'_n, Cl(W_n'))$ by:
\[F_{V_n',\varphi,1}(x\oplus h)=\]
\[
\frac{\varphi\left(\sqrt{\|x\|_1^2+\|h\|_1^2}\right)}{
\sqrt{\|x\|_1^2+\|h\|_1^2+i(h(x)-x^{*,1}(h^{*,1}))}}\left(\frac{
(h^{*,1}\oplus x^{*,1})-x\oplus h}{2}+i\frac{x\oplus h+(h^{*,1}\oplus x^{*,1})}{2}\right)
\]
for all $x\oplus h\in W_n' =V_n'\oplus (V_n')^*$.

Let $W_n '=V_n'\oplus (V_n')^*$, $W_n=V_n\otimes (V_n)^*$, and $W_n ''=V_n''\oplus (V_n'')^*$.

For each $k\geq 1$ and $l\geq 1$, we define a norm on $(\otimes ^k W_n)\otimes (\otimes ^l W_n'')$ by:

$$ \|u\|= $$
$$\sup_{\lambda_i \in W_n^*, \mu_j\in W_n'', \|\lambda_i\|\leq 1, \|\mu_j\|\leq 1, 1\leq i\leq k, 1\leq j\leq l}
((\lambda_1\otimes \cdots \otimes \lambda_k)\otimes (\mu_1\otimes \cdots \otimes \mu_l))(u)$$
for all $u\in (\otimes ^k W_n)\otimes (\otimes ^l W_n'')$.

This norm can be extended to construct Banach algebra tensor products $T(W_n)\otimes T(W_n')$ and $T_{\mathbb{C}}(W_n)\otimes T_{\mathbb{C}}(W_n')$ .
The Banach algebra norm on  $T_{\mathbb{C}}(W_n)\otimes T_{\mathbb{C}}(W_n')$ induces a Banach algebra norm on $Cl(W_n)\otimes Cl(W_n'')$
by a quotient construction.
It is not difficult to see that
$Cl(W_n)\otimes Cl(W_n')$ is naturally isomorphic  to $Cl(W_n')$.

We define a norm on the algebraic tensor product $C_0 (W_n, Cl (W_n))\otimes _{alg}
C_0(W_n '', Cl (W_n ''))$ by:

$$\| \sum_k f_k\otimes g_k\|=\sup_{w\in W_n, w''\in W_n''} \|   \sum_k (f_k (w)\otimes g_k(w''))\|$$
for all $ \sum_k f_k\otimes g_k \in C_0 (W_n, Cl (W_n))\otimes
C_0(W_n '', Cl (W_n '')). $

We define $C_0 (W_n, Cl (W_n))\otimes C_0(W_n '', Cl (W_n ''))$ to be the norm closure \linebreak
$C_0 (W_n, Cl (W_n))\otimes _{alg}
C_0(W_n '', Cl (W_n ''))$.
Observe that $C_0 (W_n, Cl (W_n))\otimes C_0(W_n '', Cl (W_n ''))$ is isomorphic
to $C_0(W_n', Cl(W_n '))$ as  Banach algebras.

For each natural number $k$,
the above Banach algebra isomorphism  can be naturally extended to a Banach algebra
 isomorphism    from $M_k(C_0 (W_n, Cl (W_n))\otimes C_0(W_n '', Cl (W_n '')))$ to
 $M_k(C_0(W_n', Cl(W_n '))).$

We can similary define the graded Banach tensor product \linebreak $C_0(W_n,Cl(W_n))\widehat{\otimes}\,C_0(W_n'',Cl(W_n''))$,
where the gradings on $C_0(W_n,Cl(W_n))$ and $C_0(W_n'',Cl(W_n''))$ are
induced by the natural gradings on the Clifford algebras.

\begin{prop} Let $F_{V_n', \varphi,1}$ be as above.
Given $r>0$ and $\epsilon>0$, there exist $\varphi$  and a natural
number $m$  such that $index(F_{V_n', \varphi,1}) \oplus p_m$ is
homotopy equivalent to $\phi_0 (D_0( (index(F_{V_n,\varphi})
\otimes index(F_{V_n'',\varphi}))\oplus (p_1\otimes p_1),
(index(F_{V_n,\varphi})\otimes p_1)\oplus (p_1\otimes
index(F_{V_n'',\varphi}))))\oplus p_m$ in $M_k ((C_0 (W_n',
Cl(W_n')))^+)$ through a homotopy of idempotents  which are
$(r,\epsilon)$-flat relative to $W_n$, where
$\phi_0$ is as in the
definition of $\phi$ in Lemma 6.6,  $D_0$ is as in the definition
of the difference construction $D$ in this section, and $p_m$ is
the direct sum of $m$ copies of the identity $1$ and $m$ copies of
$0$.
\end{prop}

{\noindent {\it Proof:}}

Let $$M= \frac{
\sqrt{\|x_1\|^2+\|h_1\|^2+i(h_1(x_1)-x_1^{*}(h_1^{*}))}
\,\,\varphi\left(\sqrt{\|x\|_1^2+\|h\|_1^2}\right)}{
\sqrt{\|x\|_1^2+\|h\|_1^2+i(h(x)-x^{*,1}(h^{*,1}))}},$$
$$N= \frac{
\sqrt{\|x_2\|^2+\|h_2\|^2+i(h_2(x_2)-x_2^{*}(h_2^{*}))}
\,\,\varphi\left(\sqrt{\|x\|_1^2+\|h\|_1^2}\right)}{
\sqrt{\|x\|_1^2+\|h\|_1^2+i(h(x)-x^{*,1}(h^{*,1}))}},$$ for all
$x=x_1\oplus x_2 \in V_n'=V_n\oplus V_n''$ and $h=h_1\oplus h_2
\in (V_n')^*=(V_n)^*\oplus (V_n'')^*.$

Let $\varepsilon$ be the grading operator in $C_b (W_n, Cl (W_n))$
induced by the natural grading of $Cl(W_n)$, where $W_n=V_n\oplus
V_n ^*$.

 Define
$$ F_{V_n', \varphi,0}= M((F_{V_n, \varphi}\otimes 1) +N( \varepsilon
\otimes F_{V_n '', \varphi})),$$ where $\varepsilon$ is the
grading operator.

Let $$F(t)=t F_{V_n', \varphi,1}+ (1-t) F_{V_n', \varphi,0}$$ for
$t\in [0,1]$.

Observe that $index (F(t))$ is a homotopy between $index (F_{V_n',
\varphi,1})$ and $index (F_{V_n', \varphi,1})$. Furthermore,
$index (F(t))$ is $(r,\epsilon)$-flat relative to $W_n$ for a
suitable choice of $\varphi$.

Let
$$M_t=\sqrt{(1-t)M^2+\frac{t}{2}},$$
$$N_t=\sqrt{(1-t)N^2+\frac{t}{2}}$$
for all $t\in [0,1]$.

It is not difficult to see that $M_t$ and $N_t$ satisfy the
conditions of Proposition 6.7. Now Proposition 6.8 follows from
Proposition 6.7 and its proof.

\begin{flushright}
$\Box$
\end{flushright}

{\noindent {\it Proof of Proposition 6.5:}} We shall prove the
$K_0$ case. The $K_1$ case can be proved in a similar way by a suspension argument.

Let $V_n'=V_n\oplus V_n''$. For any $g\in V_n^*$, we extend $g$ to
an element in $(V_n')^*$ by defining $g(x)=0$ if $x\in V_n''$.
Thus we can identify $V_n^*$ with a subspace of $(V_n')^*$.
Similarly we identify $(V_n'')^*$ with a subspace of $(V_n')^*$.
We have
$$(V_n')^*=V_n^*\oplus (V_n'')^*.$$

Define $$\|x\|_1=\sqrt{\|x_1\|^2+\|x_2\|^2}$$ for any $x=x_1\oplus
x_2\in V_n\oplus V_n''$, and
$$\|g\|_1=\sqrt{\|g_1\|^2+\|g_2\|^2}$$ for any $g=g_1\oplus g_2\in
(V_n')^*=V_n^*\oplus (V_n'')^*$.

Define
\[
\|x\|^2_t=t(\|x\|_1)^2 + (1-t)\|x\|^2\,,
\]
\[
\|g\|^2_t=t(\|g\|_1)^2 + (1-t)\|g\|^2\,,
\]
for all $t\in [0,1]$, $x=x_1\oplus x_2\in V_n\oplus V_n''$ and
$g=g_1\oplus g_2\in (V_n')^*=V_n^*\oplus (V_n'')^*$.

For each $t\in [0,1]$, let $W_n'=V_n'\oplus (V_n')^*$ be given the
norm:
\[ \|x\oplus g\|_t =\sqrt{ \|x\|_t^2 + \|g\|_t^2}\]
for all $x\oplus g \in W_n'=V_n'\oplus (V_n')^*.$

Let
\[
x^{*,t}=t(x_1^*\oplus x_2^*)+(1-t)x^*\,,
\]
\[
g^{*,t}=t(g_1^*\oplus g_2^*)+(1-t)g^*\,,
\]
for any $t\in [0,1]$, $x=x_1\oplus x_2\in V_n\oplus V_n''$,
$g=g_1\oplus g_2\in (V_n')^*=V_n^*\oplus (V_n'')^*$.

Let $\varphi$ be a continuous function on $\mathbb{R}$ such that
$0\leq \varphi(t)\leq 1$, $\exists\, 0<c_1<c_2$ satisfying
$\varphi(t)=0$ if $t\leq c_1$ and $\varphi(t)=1$ if $t\geq c_2$.

For each $t\in [0,1]$, let $F_{V_n', \varphi,t},\ \in
C_b(W_n',Cl(W_n'))$ be defined by:
\[
F_{V_n',\varphi,t}(0\oplus 0)=0,\]
\[
F_{V_n',\varphi,t}(x\oplus h)=
\]
\[
\frac{\varphi\left(\sqrt{\|x\|_t^2+\|h\|_t^2}\right)}{
\sqrt{\|x\|_t^2+\|h\|_t^2+i(h(x)-x^{*,t}(h^{*,t}))}}\left(\frac{(h^{*,t}\oplus x^{*,t})-x\oplus h}{2}
+i\frac{x\oplus h+(h^{*,t}\oplus x^{*,t})}{2}\right),
\]
for all nonzero $x\oplus h \in W_n'=V_n'\oplus (V_n')^*$.

For each $t\in [0,1]$, let $\|\cdot\|_t$ be the Banach algebra
norm on $C_b (W_n', Cl (W_n'))$ induced by the Banach space norm
$\|\cdot\|_t$ on $W_n'$.

By the definition of $index (F_{V_n', \varphi,t})$, it is not
difficult to verify that, given $r>0, \epsilon>0$, there exists
$\varphi$ such that $\varphi$ is independent of $n$, and $index (F_{V_n', \varphi,t})$ is
$(r,\epsilon)$-flat relative to $V_n\oplus 0\subseteq W_n'$ with
respect to $\|\cdot\|_t$, $\,\,\,$ i.e.

\[ \|(index (F_{V_n', \varphi,t}))(u_1)- (index (F_{V_n',
\varphi,t}))(u_2)\|_t< \epsilon\]

if $u_1, u_2\in V_n\oplus 0 \subseteq W_n'$ and $\|u_1-u_2\|_t\leq
r.$

Let
\[
index(F_{V_n',\varphi,t})=a_{V_n',\varphi, t}+\left(\begin{array}{cc} 1 & 0 \\
0 & 0
\end{array}\right)
\]
for some $a_{V_n',\varphi,t}\in C_0(W_n',Cl(W_n'))$.

Notice that
 $\|x\|_t \geq \|x\|$ for all $x\in V_n\oplus V_n''$, and
$\|g\|_t \geq \|g\|$ for all $g\in V_n^*\oplus (V_n'')^*$. This,
together with the definition of $index (F_{V_n', \varphi,t})$,
implies that there exists $R>0$ such that $R$ is independent of $n$, and

$$supp (a_{V_n',\varphi,t})\subseteq B_{W_n'}(0, R)=\{\xi\in W_n':
\|\xi\| <R\}$$ for all $n$, where $supp (a_{V_n',\varphi,t})=\{
\xi\in W_n': a_{V_n',\varphi,t}(\xi)\neq 0\}.$

For each $t\in [0,1]$, in the definition of
$\beta_{\left\{V_n'\right\}}$, we replace $F_{V_n',\varphi}$  by
$F_{V_n',\varphi,t}$ to define  $\beta_{\left\{V_n'\right\}, t}$.
Notice that
 $\|x\oplus 0\|_t= \|x\oplus 0\|$ for all $x\oplus 0\in V_n\oplus
 0 \subseteq  W_n'=V_n'\oplus (V')_n^*$, and, for each $n$, the topology
 on $C^*(F_n,\partial_\Gamma F_n,V_n)_R$ induced by $\|\cdot\|_t$ is equivalent to the
 topology induced by $\|\cdot\|.$

 The above facts, together with the
$(r,\epsilon)$-flatness of $index( F_{V_n', \varphi,t}))$ and the
above property of the support of $a_{V_n',\varphi,t}$, imply
that $\beta_{\left\{V_n'\right\}, t}$ is  a well defined
homomorphism from $K_0(C^*(\left\{F_n,\partial_\Gamma
F_n\right\}_n))$ to
$\lim_{R\to\infty}\underset{n=1}{\overset{\infty}{\oplus}} K_0(
C^*(F_n,\partial_\Gamma F_n,V_n)_R)$.

It is easy to see that
$\beta_{\left\{V_n'\right\},0}=\beta_{\left\{V_n'\right\}}$, and
$\beta_{\left\{V_n'\right\},0}=\beta_{\left\{V_n'\right\},1}$ as
homomorphisms from $K_0(C^*(\left\{F_n,\partial_\Gamma
F_n\right\}_n))$ to
$\lim_{R\to\infty}\underset{n=1}{\overset{\infty}{\oplus}} K_0(
C^*(F_n,\partial_\Gamma F_n,V_n)_R)$.

Let $W_n''=V_n''\oplus (V_n'')^*$.
Let $C_0(W_n,Cl(W_n))\widehat{\otimes}\,C_0(W_n'',Cl(W_n''))$
be  the graded Banach algebra tensor product $C_0(W_n,Cl(W_n))$ and $C_0(W_n'',Cl(W_n''))$ as defined in the paragraphs
before Proposition 6.8, where
the gradings of $C_0(W_n,Cl(W_n))$ and $C_0(W_n'',Cl(W_n''))$ are  induced by the natural gradings of the
Clifford algebras. Notice that the isomorphism from \linebreak $C_0(W_n,Cl(W_n))\widehat{\otimes}\,C_0(W_n'',Cl(W_n''))$ to
$C_0(W_n',Cl(W_n'))$ is a graded Banach algebra isomorphism.

Recall that the map: $a\widehat{\otimes} b\to a\tau^{\pt b}\otimes
b$, is an isomorphism from
\[
C_0(W_n,Cl(W_n))\widehat{\otimes}\,C_0(W_n'',Cl(W_n''))
\]
to
\[
 C_0(W_n,Cl(W_n))\otimes C_0(W_n'',Cl(W_n''))\,,
\]
where  $a$ and $b$ are respectively homogeneous elements in
$C_0(W_n,Cl(W_n))$ and $C_0(W_n'',Cl(W_n''))$,  $\tau$ is the
grading operator in $C_b(W_n,Cl(W_n))$ induced by the natural
grading of the Clifford algebra $Cl (W_n)$, and $\pt b$ is the
degree of $b$. This fact, together with Proposition 6.8, implies
that $\beta_{\{V_n'\},1}[z]$ is the product of
$\beta_{\{V_n\}}[z]$ with the direct sum of the Bott elements
associated to $\left\{W_n''\right\}_n$. Hence Proposition 6.5
follows.
\begin{flushright}
$\Box$
\end{flushright}

\section{The proof of the main result}
The purpose of  this section is to prove the main result of this
paper.

We need some preparations before we prove the main result.

Given $ d\geq 0$, let $C^*_{L,alg}(\left\{P_d(F_n)\right\}_n)$ be
the algebra of all elements
$\underset{n=1}{\overset{\infty}{\oplus}} a_n$ such that
\begin{enumerate}
\item[(1)] $a_n\in C^*_L(P_d(F_n))$\,;
\item[(2)] ${\displaystyle \sup_n}\|a_n\|<+\infty$\,;
\item[(3)] The map, $t\to \underset{n=1}{\overset{\infty}{\oplus}} a_n(t)$, is uniformly
continuous on $[0,\infty)$\,;
\item[(4)] $\,\,{\displaystyle\sup_{n,t}}\,propagation(a_n(t))<+\infty\,\,$
and $\,\,{\displaystyle\sup_n}\,propagation(a_n(t))\to 0\,\,\,$

as $t\to\infty$.
\end{enumerate}

Endow $C^*_{L,alg}(\left\{P_d(F_n)\right\}_n)$ with the norm
\[
\left\|\underset{n=1}{\overset{\infty}{\oplus}} a_n\right\|= \sup_n\|a_n\|\,.
\]

Define $C^*_{L}(\left\{P_d(F_n)\right\}_n)$ to be the norm completion
of $C^*_{L,alg}(\left\{P_d(F_n)\right\}_n)$. Let
$C^*_{L,alg}(\left\{P_d(\partial_{\Gamma}F_n)\right\}_n)$ be the
subalgebra of $C^*_{L,alg}(\left\{P_d(F_n)\right\}_n)$ consisting
of elements $\oplus_{n=1}^\infty a_n$ such that $\exists\,R\geq
0$,
\begin{eqnarray*}
supp(a_n(t))\subseteq P_d((F_n \cap B_{\Gamma}
(\Gamma-F_n,R))\times P_d( (F_n \cap B_{\Gamma} (\Gamma-F_n,R)
)\,,
\end{eqnarray*}
for any natural number $n$ and  $t\in [0,\infty)$, where
$B_{\Gamma}(\Gamma-F_n, R)=\{x\in \Gamma: d(x, \Gamma-F_n)< R\}$
if $\Gamma-F_n\neq \emptyset$, and $B_{\Gamma}(\Gamma-F_n,
R)=\emptyset$ if $\Gamma-F_n= \emptyset$.

Define $C^*_{L}(\left\{ P_d(\pt_\Gamma F_n)\right\}_n)$ to be the
norm closure of $C^*_{L,alg}(\left\{ P_d(\pt_\Gamma
F_n)\right\}_n)$. Note that $C^*_{L}(\left\{ P_d(\pt_\Gamma
F_n)\right\}_n)$ is a two-sided ideal of
$C^*_{L}(\left\{P_d(F_n)\right\}_n)$.

Throughout the rest of this paper, we fix $x_0\in \Gamma$ and choose
$$F_n=B(x_0, n)=\{x\in \Gamma: d(x,x_0)\leq n\}$$
for each non-negative integer $n$.

For any $d\geq 0$, let $\partial_{\Gamma, d} F_n$ be
$$\{x\in \Gamma: n-10d \leq d(x,x_0)\leq n\}.$$

By the choice of $F_n$,  $P_d (F_n)\backslash P_d(\partial_{\Gamma,d} F_n)$
is an open subset of $P_d (F_{n+1})\backslash P_d(\partial_{\Gamma,d} F_{n+1})$.
Let $i_{n,d}$ be the inclusion homomorphism from $ C_0(P_d (F_n)\backslash P_d(\partial_{\Gamma,d} F_n))$ into
$C_0(P_d (F_{n+1})\backslash P_d(\partial_{\Gamma,d} F_{n+1}))$.
By the definition of relative K-homology group in Section 4, $i_{n,d}$ induces  a  homomorphism
$$(i_{n,d})^*:\,\,
K_* (P_d (F_{n+1}), P_d(\partial_{\Gamma,d} F_{n+1})) \rightarrow
K_* (P_d (F_n), P_d(\partial_{\Gamma,d} F_n)).$$

We have $$K_*(\prod_{n} P_d (F_n),\prod_{n}P_d(\partial_{\Gamma,d} F_n))=\oplus_{n}K_* (P_d (F_n), P_d(\partial_{\Gamma,d} F_n)),$$
where $\oplus_{n}K_* (P_d (F_n), P_d(\partial_{\Gamma,d} F_n))$ is defined to be
$$\{\oplus_n z_n:\,\, z_n \in K_* (P_d (F_n), P_d(\partial_{\Gamma,d} F_n))\}.$$

For each $d\geq 0$, we define a homomorphism
$$s_d: \,\,K_*(\prod_{n} P_d (F_n),\prod_{n}P_d(\partial_{\Gamma,d} F_n))\rightarrow
K_*(\prod_{n} P_d (F_n),\prod_{n}P_d(\partial_{\Gamma,d} F_n))$$ by:

$$ s_d(\oplus_n z_n) =  \oplus_{n} (i_{n,d})^* (z_{n+1})$$
for all $\oplus_n z_n \in K_*(\prod_{n} P_d (F_n),\prod_{n}P_d(\partial_{\Gamma,d} F_n))$.

Notice that if $d_1>d_2\geq 0$, then
$$s_{d_1}(\oplus_n z_n)=s_{d_2}(\oplus_n z_n) \in
\lim_{d\rightarrow \infty} K_*(\prod_{n} P_d (F_n),\prod_{n}P_d(\partial_{\Gamma,d} F_n))$$
for all $\oplus_n z_n \in K_*(\prod_{n} P_{d_2} (F_n),\prod_{n}P_{d_2}(\partial_{\Gamma,d_2} F_n))$.

It follows that $\{s_d\}_d$ induces a well defined homomorphism
$$s:\,\,
\lim_{d\rightarrow \infty} K_*(\prod_{n} P_d (F_n),\prod_{n}P_d(\partial_{\Gamma,d} F_n))\rightarrow
\lim_{d\rightarrow \infty} K_*(\prod_{n} P_d (F_n),\prod_{n}P_d(\partial_{\Gamma,d} F_n)).$$

Observe that $P_d (F_n)\backslash P_d(\partial_{\Gamma,d} F_n)$ is an open subset of $P_d(\Gamma)$.
Let $j_{n,d}$ be the inclusion homomorphism from $C_0( P_d (F_n)\backslash P_d(\partial_{\Gamma,d} F_n))$
to $C_0(P_d (\Gamma))$. $j_{n,d}$ induces a homomorphism
$$(j_{n,d})^*:\,\, K_*(P_d(\Gamma))\rightarrow
K_*(P_d (F_n),P_d(\partial_{\Gamma,d} F_n)).$$

 We define a homomorphism
$$r: \lim_{d\rightarrow \infty}K_*(P_d (\Gamma))\rightarrow \lim_{d\rightarrow \infty} K_*(\prod_{n} P_d (F_n),
\prod_{n}P_d (\partial_{\Gamma,d} F_n)) $$ by:
$$ r( z)= \oplus_n (j_{n,d})^* z$$
for all $z\in K_*(P_d (\Gamma)).$

\begin{prop}
We have the following exact sequence:
%\begin{eqnarray*}
$$ \rightarrow  \lim_{d\rightarrow \infty}K_{*+1}(\prod_{n} P_d (F_n),\prod_{n}P_d ( \partial_{\Gamma,d}
F_n))\overset{Id-s}{\rightarrow} \lim_{d\rightarrow \infty}K_{*+1}(\prod_{n} P_d
(F_n),\prod_{n}P_d (
\partial_{\Gamma,d} F_n)) \rightarrow $$ $$\lim_{d\rightarrow \infty}K_*(P_d (\Gamma))
 \overset{r }{\rightarrow} \lim_{d\rightarrow \infty} K_*(\prod_{n} P_d (F_n),
\prod_{n}P_d (\partial_{\Gamma,d} F_n)) \overset{Id-s}{\rightarrow}
\lim_{d\rightarrow \infty} K_*(\prod_{n} P_d (F_n),
\prod_{n}P_d(\partial_{\Gamma,d} F_n)) \rightarrow ,$$
where $Id$ is the identity homomorphism.
%\end{eqnarray*}
\end{prop}

The above proposition follows from the standard $lim^1$-sequence
for $K$-homology [36].

\begin{defi}
We define the $C^*$-algebra $C^*_{L,d}(\left\{F_n,\pt_\Gamma
F_n\right\}_n)$ to be the quotient algebra of
$C^*_{L}(\left\{P_d(F_n)\right\}_n)$ over $C^*_{L}(\left\{
P_d(\pt_\Gamma F_n)\right\}_n)$.
\end{defi}

We   define a map
\[
\chi_L\,:\,C^*_L(P_d (\Gamma))\to C^*_{L,d} (\left\{F_n,
\pt_\Gamma F_n\right\})
\]
by
\[
\chi_L (a)=[\underset{n=1}{\overset{\infty}{\oplus}}\chi_n a
\chi_n]\,
\]
for all $a\in C^*(\Gamma),$ where $\chi_n$ is the characteristic
function of $F_n$.

\begin{prop}
We have the following exact sequence:
%\begin{eqnarray*}
$$\rightarrow  \lim_{d\rightarrow \infty}K_{*+1}(C^*_{L,d}(\left\{F_n,\pt_\Gamma
F_n\right\}_n))\overset{(Id-S)_*}{\rightarrow} \lim_{d\rightarrow
\infty} K_{*+1}(C^*_{L,d}(\left\{F_n,\pt_\Gamma F_n\right\}_n))
\rightarrow
$$ $$\lim_{d\rightarrow\infty}K_*(C^*_{L}(P_d(\Gamma)))
 \overset{ (\chi_L)_*}{\rightarrow} \lim_{d\rightarrow \infty}K_*(C^*_{L,d}(\left\{F_n,\pt_\Gamma
F_n\right\}_n)) \overset{(Id-S)_*}{\rightarrow}
\lim_{d\rightarrow\infty}K_*(C^*(\left\{F_n,\pt_{\Gamma} F_n\right\}_n))
\cdots,$$ where $Id$ is the identity homomorphism.
%\end{eqnarray*}
\end{prop}

{\noindent {\it Proof:}} We have the following commutative
diagram:

\vspace{0.50cm}

\begin{flushleft}
\begin{tabular}{cccc}
${\displaystyle \lim_{d\rightarrow \infty}K_{*+1}(\prod_{n} P_d
(F_n),\prod_{n} P_d ( \partial_{\Gamma,d} F_n))}$ &
${\displaystyle\overset{Id-s}{\rightarrow}}$ &
${\displaystyle\lim_{d\to\infty}K_{*+1}(\prod_{n} P_d
(F_n),\prod_{n} P_d (
\partial_{\Gamma,d} F_n) )}$
&${\displaystyle\overset{}{\rightarrow}}$   \\
${\displaystyle\downarrow }$ & & ${\displaystyle\downarrow }$
\\ ${\displaystyle \lim_{d\rightarrow \infty}K_{*+1}(C^*_{L,d}(\left\{F_n,\pt_\Gamma
F_n\right\}))}$ & ${\displaystyle
\overset{(Id-S)_*}{\rightarrow}}$ & ${\displaystyle
\lim_{d\rightarrow \infty}K_{*+1}(C^*_{L,d}(\left\{F_n,\pt_\Gamma
F_n\right\}))}$
&${\displaystyle\overset{}{\rightarrow}}$ \\
 & &
\\ & &
&
\end{tabular}
\end{flushleft}

\begin{tabular}{ccccc}
${\displaystyle\overset{}{\rightarrow}}$ &
${\displaystyle\lim_{d\to\infty}K_*(P_d(\Gamma))}$ &
${\displaystyle\overset{r}{\rightarrow}}$ &
${\displaystyle\lim_{d\to\infty}K_*(\prod_{n} P_d
(F_n),\prod_{n}P_d (
\partial_{\Gamma,d} F_n) )}$& ${\displaystyle\overset{Id-s}{\rightarrow}}$\\  &
${\displaystyle\downarrow }$ & &
${\displaystyle \downarrow }$ & \\
${\displaystyle\overset{}{\rightarrow}}$ & ${\displaystyle
\lim_{d\rightarrow \infty}K_*(C^*_L(P_d (\Gamma)))}$ &
${\displaystyle\overset{(\chi_L)_*}{\rightarrow}}$ &
${\displaystyle \lim_{d\rightarrow
\infty}K_*(C^*_{L,d}(\left\{F_n,\pt_\Gamma F_n \right\}))}$&
${\displaystyle \overset{ (Id-S)_*}{\rightarrow}}$\\
& & & \\
&  &
\end{tabular}

By Theorem 4.5 and Proposition 4.6, the vertical maps in the above
diagram are isomorphisms. Notice that the first horizontal
sequence is exact. This, together with the commutativity of the
diagram, implies that the second horizontal sequence is exact.

\begin{flushright} $\Box$
\end{flushright}

Notice that $\lim_{d\rightarrow \infty}
C^*(\left\{P_d(F_n)\right\}_n)$ is naturally $*$-isomorphic to
$C^*(\left\{F_n\right\}_n)$ (cf. [25]).

 Let $e$ be the evaluation homomorphism from
$C^*_{L}(\left\{P_d(F_n)\right\}_n)$ to
$C^*(\left\{F_n\right\}_n)$ defined by:
$$e(\underset{n=1}{\overset{\infty}{\oplus}}
a_n)=\underset{n=1}{\overset{\infty}{\oplus}} a_n(0)$$ for all
$\underset{n=1}{\overset{\infty}{\oplus}}a_n \in
C^*_{L}(\left\{P_d(F_n)\right\}_n)$, where
$\underset{n=1}{\overset{\infty}{\oplus}} a_n(0)$ is identified
with an element of $C^*(\left\{F_n\right\}_n)$ by the above
isomorphism.

$e$ induces a $\ast$-homomorphism (still denoted by $e$) from
$C^*_{L,d}(\left\{F_n,\pt_\Gamma F_n\right\}_n)$ to
$C^*(\left\{F_n,\pt_\Gamma F_n\right\}_n)$.

\begin{lem}
$\beta_{\left\{V_n\right\}}\circ e_*$ is an isomorphism of
$\lim_{d\to\infty} K_*(C^*_{L}(\left\{P_d(F_n)\right\}_n))$ onto
$\lim_{R\to\infty}\underset{n=1}{\overset{\infty}{\oplus}}
K_*(C^*(F_n,V_n)_R)$.
\end{lem}

We need some preparations to prove Lemma 7.4.

For any $d\geq 0$, let $C^* (P_d(F_n))\otimes_{alg} C_0(W_n, Cl(W_n))$ be the algebraic tensor product of
$C^* (P_d(F_n))$ with $ C_0(W_n, Cl(W_n))$. We can construct a representation of
the algebra $C^* (P_d(F_n))\otimes_{alg} C_0(W_n, Cl(W_n))$ on the Banach space $C_0 (W_n, (l^2(\Gamma_{d,n})\otimes H) \otimes Cl(W_n)),$
where $C^*(P_d (F_n))$ is defined as in Section 4 using  a
countable dense subset $\Gamma_{d,n}$ of $P_d (F_n)$, and $C_0 (W_n, (l^2(\Gamma_{d,n})\otimes H) \otimes Cl(W_n))$ is the Banach
space defined as in the definition of $C^*(F_n, V_n)$.
We define $C^* (P_d(F_n))\otimes C_0(W_n, Cl(W_n))$ to be the operator norm closure of $C^* (P_d(F_n))\otimes_{alg} C_0(W_n, Cl(W_n))$,
where the operator norm is induced by the representation of
the algebra $C^* (P_d(F_n))\otimes_{alg} C_0(W_n, Cl(W_n))$ on the Banach space $C_0 (W_n, (l^2(\Gamma_{d,n})\otimes H) \otimes Cl(W_n))$.

For any
$\eta$ in $l^{\infty} (\Gamma_{d,n})$, we can define an operator $N_{\eta}$ on
the Banach space $C_0 (W_n, (l^2(\Gamma_{d,n})\otimes H) \otimes Cl(W_n))$
by:
$$(N_{\eta}\zeta)(w)= ((\eta\otimes I)\otimes 1) (\zeta(w) ) $$
for all $\zeta \in C_0 (W_n, (l^2(\Gamma_{d,n})\otimes H) \otimes Cl(W_n))$ and $w\in W_n$,
where $\eta$ acts on $l^2(\Gamma_{d,n})$ by multiplication.

We can verify that $N_{\eta}$ is a bounded
operator (using  an argument similar to the proof of Claim 2 in the proof of Lemma 6.1).

For every pair $(x,y)$ in $\Gamma_{d,n}
\times \Gamma_{d,n}$, let $\delta_x$ and $\delta_y$ be
respectively the Dirac functions in $l^{\infty}(\Gamma_{d,n})$ at
$x$ and $y$. For any $a\in C^*(P_d (F_n), V_n)$, we write
$$ a(x,y)= N_{\delta_x} a N_{\delta_y} .$$

We  define
$$ supp (a)=\{ (x,y)\in \Gamma_d \times \Gamma_d \,|\, a(x,y)\neq
0\},$$
$$propagation (a)= sup\{ d(x,y)\,|\, (x,y)\in supp (a)\}.$$

For any
 $g$ in $C_0(W_n)$ and $a$ in $C^*(P_d (F_n), V_n)  $,
we define a bounded operator $ga$ on the Banach space $C_0 (W_n, (l^2(\Gamma_{d,n})\otimes H) \otimes Cl(W_n))$
by:
$$((g a)\zeta)(w)= g(w)\,\,( a \zeta)(w)  $$
for all $\zeta \in C_0 (W_n, (l^2(\Gamma_{d,n})\otimes H) \otimes Cl(W_n))$ and $w\in W_n$.

We define the support of $a(x,y)$, $support (a(x,y))$, to be the complement  of the set of
all points $\xi \in W_n$ such that there exists $g\in C_0(W_n)$ satisfying $g(\xi)\neq 0$
and $(ga)(x,y)=0$.

We define $C^*_{alg, L}(P_d (F_n), V_n) $ to be the algebra of all
uniformly continuous functions
$$a: [0,\infty) \rightarrow
C^*(P_d (F_n), V_n)$$ satisfying
$\,\,\,\,{\displaystyle\sup_{t}}\,propagation(a(t))<+\infty\,\,\,$ and
$\,\,\, propagation(a(t))\to 0\,\,\,\,$
as $t$ goes to $\infty$.

We endow  $C^*_{alg, L} (P_d(F_n), V_n)$ with the norm:
$$\| a\| = sup_{t\in [0,\infty)} \| a(t)\|.$$

We define $C^*_{ L} (P_d(F_n), V_n)$  to be the norm closure of
$C^*_{alg, L} (P_d(F_n), V_n)$.

For each $n$ and $R>0$, we define $C^*_{ L} (P_d(F_n), V_n)_R$ to
be the  closed subalgebra of $C^*_L(P_d (F_n), V_n)$ generated by all
elements $a\in C^*_L(P_d (F_n), V_n)$ such that
\begin{enumerate}

\item[(1)]
 $propagation
(a(t))<R$ for all $t\in [0,\infty)$

\item[(2)] $support
((a(t))(x,y))\subseteq B_{W_n}(f(x)\oplus 0, R)$ for all $t\in
[0,\infty)$, all $x$ and $y$ in $\Gamma_{d,n}$, where $f: P_d (\Gamma) \rightarrow X,$ is the convex linear
extension of the uniform embedding $f: \Gamma \rightarrow  X$.
\end{enumerate}

We can define the Bott map
\[
\beta_{\left\{V_n\right\},L}\,:\,K_*(C^*_L(\left\{P_d(F_n)\right\}_n))\to
\lim_{R\rightarrow\infty}\underset{n=1}{\overset{\infty}{\oplus}}
K_*(C^*_L(P_d(F_n),V_n)_R)\,,
\]
in a way similar to the definition of
$\beta_{\left\{V_n\right\}}$.
\begin{lem}
$\,\,\,\,\beta_{\left\{V_n\right\},L}\,\,$ is an isomorphism between
$\,\,\,K_*(C^*_L(\left\{P_d(F_n)\right\}_n))\,\,\,\,\,$ and
$\,\,\,\,\,\,\,\,\lim_{R\rightarrow\infty}\underset{n=1}{\overset{\infty}{\oplus}}
K_*(C^*_L(P_d(F_n),V_n)_R) \,. $
\end{lem}
{\noindent {\it Proof:}} By a Mayer-Vietoris sequence argument and
induction on the dimension of the skeletons, the general case can
be reduced to the 0-dimensional case, i.e., if $\Gamma_n\subseteq
P_d(F_n)$ is $\delta$-separated (meaning $d(x,y)\geq \delta$ if
$x\not= y\in \Gamma_n$) for some $\delta>0$, then
$\beta_{\left\{V_n\right\},L}$ is an isomorphism from
$K_*(C^*_{L}(\left\{\Gamma_n\right\}_n))$ to
$K_*(C^*_{L}(\left\{\Gamma_n,V_n\right\}))$\,. This 0-dimensional case follows
from the facts that
\[
K_*(C^*_{L}(\left\{\Gamma_n\right\}_n))\cong
\underset{n=1}{\overset{\infty}{\oplus}}\underset{z\in\Gamma_n}{\oplus}K_*(C^*_L(\left\{z\right\}))\,,
\]
\[
\lim_{R\rightarrow\infty}\underset{n=1}{\overset{\infty}{\oplus}}K_*(C^*_{L}(\Gamma_n,V_n)_R)\cong
\underset{n=1}{\overset{\infty}{\oplus}}\underset{z\in\Gamma_n}{\oplus}K_*(C^*_L(\left\{z\right\},V_n))\,,
\]
and $\beta_{\left\{V_n\right\},L}$ is an isomorphism from
$K_*(C^*_L(\left\{z\right\}))$ to
$K_*(C^*_L(\left\{z\right\},V_n))$ by Lemma 3.1 and the Bott
periodicity.
\begin{flushright}
$\Box$
\end{flushright}

{\noindent {\it Proof of Lemma 7.4:}}

Notice that $\lim_{d\rightarrow \infty}C^*(P_d(F_n),V_n)_R $ is
naturally isomorphic to $C^*(F_n,V_n)_R$.

Let $e_{V_n}$ be the homomorphism:
\[
C^*_L(P_d(F_n),V_n)_R \to  C^*(F_n,V_n)_R
\]
defined by
\[
e_{V_n}( a)= a(0)
\]
for each $ a \in C^*_L(P_d(F_n),V_n)_R$, where $a(0)$ is
identified with an element in $C^*(F_n,V_n)_R$ by the above
isomorphism.

The family of homomorphisms $\{e_{V_n}\}$ induces a homomorphism
$\left(e_{\left\{V_n\right\}}\right)_*$:
\[
\lim_{d\to\infty}\lim_{R\to \infty}
\underset{n=1}{\overset{\infty}{\oplus}}
K_*(C^*_L(P_d(F_n),V_n)_R)\to \lim_{R\to \infty}
\underset{n=1}{\overset{\infty}{\oplus}} K_*(C^*(F_n,V_n)_R)\,.
\]
Clearly we have
\[
\left(e_{\left\{V_n\right\}}\right)_*\circ
\beta_{\left\{V_n\right\},L}= \beta_{\left\{V_n\right\}}\circ e_*\
.
\]
By Lemma 7.5, it is enough to prove that
$\left(e_{\left\{V_n\right\}}\right)_*$ is an isomorphism. Assume
that for each $n$, $O_n$ is an open subset of $W_n=V_n\oplus
V_n^*$. We define $C^*(F_n,O_n)$ to be a closed subalgebra of
$C^*(F_n,V_n)$ consisting of elements $a\in C^*(F_n,V_n)$ such
that $support(a(x,y))\subseteq O_n$ for all $x,y\in F_n$.

Similarly we can define $C^*(F_n,O_n)_R$ for each $R>0$. For each
$r>0$, let
\[
O_{n,r}= \bigcup_{x\in F_n} B_{W_n}(f(x)\oplus 0\,,\, r)\,.
\]
By the bounded geometry property of $\Gamma$ and the fact that $f$
is a uniform embedding, there exists a natural number $m$
(independent of $n$) such that
\begin{enumerate}
\item[(1)] $O_{n,r}={\displaystyle \bigcup_{k=1}^m}
O_{n,r}^{(k)}$\,;
\item[(2)] for each $k$, $O_{n,r}^{(k)}$ is a disjoint union of
open balls of the form $B_{W_n}(f(x)\oplus 0,r)$.
\end{enumerate}
By a Mayer-Vietoris sequence argument, it suffices to prove that
$(e_{\left\{V_n\right\}})_*$ is an isomorphism:
\[
\lim_{d\to\infty}\lim_{R\to\infty}\underset{n=1}{\overset{\infty}{\oplus}}
K_*(C^*_{L}(P_d(F_n),O_{n,r}^{(k)})_R)
\to\lim_{R\to\infty}\underset{n=1}{\overset{\infty}{\oplus}}K_*(C^*(F_n,O_{n,r}^{(k)})_R)
\]
for each $k$ and $r$.

Assume that $O_{n,r}^{(k)}$ is the disjoint union of
$B_{W_n}(f(z)\oplus 0,r)$ ($z\in F_{n,r}^{(k)}$), where
$F_{n,r}^{(k)}$ is a subset of $F_n$.

For each $D>0$, we define $A_{D,n,r}^{(k)}$ to be the Banach
algebra consisting of elements $\oplus_{z\in F_{n,r}^{(k)}}\, b_z$
such that
\begin{enumerate}
\item[(1)] $b_z\in C^*(B_{F_n}(z,D),V_n)$, where $B_{F_n}(z,D)=\left\{z\in F_n\,:\,
d(x,z)<D\right\}$;
\item[(2)] $support(b_z (x,y))\subseteq B_{W_n}(f(z)\oplus 0,r)$
for all $x,\,y \in B_{F_n}(z,D)$.
\end{enumerate}

$A_{D,n,r}^{(k)}$ is endowed with the norm
\[
\left\|\underset{z\in F_{n,r}^{(k)}}{\oplus}b_z\right\|=\sup_{z\in
F_{n,r}^{(k)}}\|b_z\|\,.
\]

It is easy to see that
\[
\lim_{R\to\infty}\underset{n=1}{\overset{\infty}{\oplus}}K_*(C^*(F_n,O_{n,r}^{(k)})_R)
\]
is naturally isomorphic to
\[
\lim_{D\to
\infty}\underset{n=1}{\overset{\infty}{\oplus}}K_*\left(A_{D,n,r}^{(k)}\right)\,.
\]

For each $n$,  we  define $A_{L,D,n,r}^{(k)}$ to be the Banach
algebra generated by elements $\underset{z\in
F_{n,r}^{(k)}}{\oplus} a_{z}$ such that
\begin{enumerate}
\item[(1)] $a_{z}\in C^*_{L}(P_{D'}(B_{F_n}(z,D)),V_n)$\, where $D' =diameter(B_{F_n} (z, D))$;
\item[(2)] ${\displaystyle \sup_{z \in F_{n,r}^{(k)}}\|a_{z}\|<+\infty}$\,;
\item[(3)] $support(a_{z}(x,y))\subseteq B_{W_n}(f(z)\oplus 0,r)$ for all
$z \in F_{n,r}^{(k)}$ and $x,\,y \in B_{F_n}(z,D)$\,;
\item[(4)] $t\to \underset{z\in
F_{n,r}^{(k)}}{\oplus} a_{z}(t)$ is uniformly continuous;
\item[(5)] ${\displaystyle
\sup_{z}\,propagation(a_{z}(t))<+\infty}$, and ${\displaystyle
\sup_{z}\,propagation(a_{z}(t))\to 0}$ as \\$t\to\infty$, where
$A_{L,D,n, r}^{(k)}$ is endowed with the norm
\[
\left\|\underset{z\in F_{n,r}^{(k)}}{\oplus}
a_{z}\right\|=\sup_{z\in F_{n,r}^{(k)}}\|a_{z}\|\,.
\]

It is easy to see that
\[
\lim_{d\to\infty}\lim_{R\rightarrow\infty}\underset{n=1}{\overset{\infty}{\oplus}}
K_*(C^*_L(P_d(F_n),O_{n,r}^{(k)})_R)
\]
is naturally isomorphic to
\[
\lim_{D\to
\infty}\underset{n=1}{\overset{\infty}{\oplus}}K_*\left(A_{L,D,n,r}^{(k)}\right)\,.
\]
Using the fact that $P_{D'}( B_{F_n}(z, D))$ is (Lipschitz)
contractible,  it is easy to see that the evaluation homomorphism
(at 0) induces an isomorphism from
\[\lim_{D\to
\infty}\underset{n=1}{\overset{\infty}{\oplus}}K_*\left(A_{L,D,n,r}^{(k)}\right)
\]
to
\[\lim_{D\to
\infty}\underset{n=1}{\overset{\infty}{\oplus}}K_*\left(A_{D,n,r}^{(k)}\right)\,.
\]
This implies Lemma 7.4.
\end{enumerate}
\begin{flushright} $\Box$
\end{flushright}

The following result can be proved in a way similar to Lemma 7.4.
\begin{lem}
$\beta_{\{V_n\}}\circ e_*$ is an isomorphism between
\[
\lim_{d\to \infty}K_*\left(C^*_L(\left\{ P_d(\pt_\Gamma
F_n)\right\}_n)\right)
\]
and
\[
\lim_{R\to
\infty}\underset{n=1}{\overset{\infty}{\oplus}}K_*\left(C^*\left(\{\pt_\Gamma
F_n,V_n\} \right)_R\right)\,.
\]
\end{lem}

\begin{prop}
$\beta_{\{V_n\}}\circ e_*$ is an isomorphism between
\[
\lim_{d\to \infty}K_*\left(C^*_{L,d}(\{F_n,\pt_\Gamma
F_n\}_n)\right)
\]
and
\[
\lim_{R\to
\infty}\underset{n=1}{\overset{\infty}{\oplus}}K_*\left(C^*\left(F_n,\pt_\Gamma
F_n,V_n \right)_R\right)\,.
\]
\end{prop}
{\noindent {\it Proof:}} The proposition follows from Lemmas 7.4
and 7.6,  and a five lemma argument.
\begin{flushright} $\Box$
\end{flushright}
{\noindent {\it Proof of Theorem 1.1:}}

Let $\chi_L$ be as in Proposition 7.3. Let $e$ be the evaluation
homomorphism from $C^*_L (P_d(\Gamma))$ to $C^*(P_d (\Gamma))$
defined by: $e(a)=a(0).$ $e$ induces a homomorphism (still denoted
by $e_*$)  from $K_*(C^*_L (P_d(\Gamma)))$ to $K_*(C^* (\Gamma)).$

Let $e'_*= \chi_* \circ e_*$ be the homomorphism from $K_*(C^*_L
(P_d(\Gamma)))$ to  $K_*(C^*_S(\Gamma)),$  where $\chi$ is as in
Lemma 5.6. It is enough to prove that $e'_*$ is injective from
$\lim_{d\to\infty}K_* (P_d(\Gamma))$ to $K_* (C^*_S (\Gamma)).$

 We have the
following commutative diagram:

\vspace{0.50cm}

\begin{flushleft}
\begin{tabular}{cccc}
${\displaystyle \lim_{d\to\infty }
K_{*+1}(C^*_{L,d}(\left\{F_n,\pt_\Gamma F_n\right\}))}$ &
${\displaystyle\overset{(Id -S)_*}{\longrightarrow}}$ &
${\displaystyle\lim_{d\to\infty}K_{*+1}(C^*_{L,d}(\left\{F_n,\pt_\Gamma
F_n\right\}))}$
&${\displaystyle\overset{\gamma}{\longrightarrow}}$   \\
${\displaystyle\downarrow e_*}$ & & ${\displaystyle\downarrow
e_*}$ \\ ${\displaystyle K_{*+1}(C^*(\left\{F_n,\pt_\Gamma
F_n\right\}))}$ & ${\displaystyle
\overset{(Id-S)_*}{\longrightarrow}}$ & ${\displaystyle
K_{*+1}(C^*(\left\{F_n,\pt_\Gamma F_n\right\}))}$
&${\displaystyle\overset{\gamma}{\longrightarrow}}$ \\
${\displaystyle\downarrow \beta_{\left\{V_n\right\}}}$ & &
${\displaystyle\downarrow \beta_{\left\{V_n\right\}}}$\\
${\displaystyle\lim_{R\to\infty }
\underset{n=1}{\overset{\infty}{\oplus}}
K_{*+1}(C^*(F_n,\pt_\Gamma F_n,V_n)_R)}$ & &
${\displaystyle\lim_{R\to\infty}\underset{n=1}{\overset{\infty}{\oplus}}K_{*+1}(C^*(F_n,\pt_\Gamma
F_n,V_n)_R)}$
\end{tabular}
\end{flushleft}

\vspace{1.25cm}

\begin{tabular}{cccc}
${\displaystyle\overset{\gamma}{\longrightarrow}}$ &
${\displaystyle\lim_{d\to\infty}K_*(C^*_{L}(P_d(\Gamma)))}$ &
${\displaystyle\overset{(\chi_L)_*}{\longrightarrow}}$ &
${\displaystyle\lim_{d\to\infty}K_*(C^*_{L,d}(\left\{F_n,\pt_\Gamma
F_n \right\}))}$\\  & ${\displaystyle\downarrow e'_*}$ & &
${\displaystyle \downarrow e_*}$\\
${\displaystyle\overset{\gamma}{\longrightarrow}}$ &
${\displaystyle K_*(C^*_S(\Gamma))}$ &
${\displaystyle\overset{j_*}{\longrightarrow}}$ &
${\displaystyle K_*(C^*(\left\{F_n,\pt_\Gamma F_n \right\}))}$\\
& & & ${\displaystyle\downarrow \beta_{\left\{V_n\right\}}}$\\ & &
&
${\displaystyle\lim_{R\to\infty}\underset{n=1}{\overset{\infty}{\oplus}}K_*(C^*(F_n,\pt_\Gamma
F_n,V_n)_R)}$
\end{tabular}

\vspace{1.00cm}

By Propositions 7.3 and 5.5,  the first and second horizontal
sequences in the above diagram are exact.

Let $[x]\in \lim_{d\to\infty}K_*(C^*_L(P_d(\Gamma)))$ such that
$e'_*[x]=0$. We need to prove that $[x]=0$.

We first claim that $(\chi_L)_*[x]=0$. This follows from the
identity
\[
\beta_{\{V_n\}}(e_*((\chi_L)_*[x]))=\beta_{\{V_n\}}(j_*(e'_*[x]))=0
\]
(by the commutativity of the diagram), and the fact that
$\beta_{\{V_n\}}\circ e_*$ is an isomorphism (by Proposition 7.7).

By exactness, $\exists\,[x']\in
\lim_{d\to\infty}K_{*+1}(C^*_{L,d}(\{F_n,\pt_\Gamma F_n\}))$ such
that $\gamma[x']=[x]$. By the commutativity of the diagram, we
have
\[
\gamma (e_*[x'])=e'_*(\gamma[x'])=e'_*[x]=0\ .
\]
Hence by exactness, $\exists\,[y]\in K_{*+1}(C^*(\{F_n,\pt_\Gamma
F_n\}))$ such that $$(Id-S)_*[y]=e_*[x'].$$ The fact that
$\beta_{\{V_n\}} \circ e_*$ is an isomorphism implies that
\[
\exists\,[x'']\in \lim_{d\to\infty}K_{*+1}(C^*_{L,d}(\{F_n,\pt_\Gamma F_n\}))
\]
such that $$\beta_{\{V_n\}}\circ e_*[x'']=\beta_{\{V_n\}}[y].$$
Hence we have
\[
\beta_{\{V_n\}}( e_*[x'']-[y])=0\,.
\]
This, together with Proposition 6.5, implies that
\[
\beta_{\{V_n\}}(Id-S)_*( e_*[x'']-[y])=0\,.
\]
By the commutativity of the diagram and the property of $y$, it
follows that
\[
\beta_{\{V_n\}}(e_*(Id-S)_*[x'']-e_*[x'])=0\,.
\]
Hence
\[
\beta_{\{V_n\}}\circ e_*( (Id-S)_*[x'']-[x'])=0\,.
\]
By Proposition 7.7, we have
\[(Id-S)_*[x'']-[x']=0\,.
\]
This,  together with exactness and the identity $[x]= \gamma[x']$,
implies that $[x]=0$.
\begin{flushright} $\Box$
\end{flushright}

\noindent Department of Mathematics, 1326 Stevenson Center,

\noindent Vanderbilt University, Nashville, TN 37240, USA

\noindent e-mail: gyu@math.vanderbilt.edu and
kasparov@math.vanderbilt.edu


\begin{thebibliography}{99}

\bibitem { } E. Asplund, Averaged norms. Israel J. Math. 5 (1967), 227--233.


\bibitem { } P. Baum and A. Connes, K-theory for discrete groups,
Operator Algebras and Applications, (D. Evans and M. Takesaki,
editors), Cambridge University Press (1989), 1--20.

\bibitem { } N. Brown and E. Guentner, Uniform embedding of
bounded geometry spaces into  reflexive Banach spaces, to appear in Proc. AMS.



\bibitem { } A. Connes, Cyclic cohomology and transverse
fundamental class of a foliation, Geometric Methods in Operator
Algebras, (H. Araki and E.G. Effros, editors), Pitman Res. Notes
Math., Vol. 123 (1996), 52--144.

\bibitem { } A. Connes, Noncommutative Geometry, Academic Press,
1994.

\bibitem { } A. Connes, M.  Gromov and H. Moscovici, Conjecture de
Novikov et fibrés presque plats. C. R. Acad. Sci. Paris Sér. I
Math. 310 (1990), no. 5, 273--277.

\bibitem { } A. Connes, M. Gromov and H. Moscovici, Group
cohomology with Lipschitz control and higher signatures, Geometric
and Functional Analysis, 3 (1993), 1--78.

\bibitem { } A. Connes and H. Moscovici, Cyclic cohomology, the
Novikov conjecture and hyperbolic groups, Topology 29 (1990),
345--388.

\bibitem { } A. Connes and G. Skandalis,  The longitudinal index theorem
for foliations. Publ. Res. Inst. Math. Sci. 20 (1984), no. 6, 1139--1183.



\bibitem { } J. Diestel,  Geometry of Banach spaces---selected topics.
Lecture Notes in Mathematics, Vol. 485. Springer-Verlag,
Berlin-New York, 1975.

\bibitem { } A. N. Dranishnikov, S. C. Ferry, S. Weinberger,  Large
Riemannian manifolds which are flexible. Ann. of Math. (2) 157
(2003), no. 3, 919--938.

\bibitem { } A. N. Dranishnikov, G. Gong, V. Lafforgue, and G. Yu,  Uniform
embeddings into Hilbert space and a question of Gromov. Canad.
Math. Bull. 45 (2002), no. 1, 60--70.

\bibitem { } P. Enflo,  Banach spaces which can be given an equivalent uniformly convex norm,
Israel J. Math. 13 (1972), 281--288.




\bibitem { } G. Gong and G. Yu,  Volume growth and positive scalar
curvature. Geom. Funct. Anal. 10 (2000), no. 4, 821--828.

\bibitem { } M. Gromov, Asymptotic invariants for infinite
groups, Geometric Group Theory, (G. A. Niblo and M. A. Roller,
editors), Cambridge University Press, (1993), 1--295.

\bibitem { } M. Gromov, Problems (4) and (5), Novikov
Conjectures, Index Theorems and Rigidity, Vol. 1, (S. Ferry, A.
Ranicki and J. Rosenberg, editors), Cambridge University Press,
(1995), 67.

\bibitem { } M. Gromov, Positive curvature, macroscopic
dimension, spectral gaps and higher signatures, Functional
Analysis on the eve of the 21st century, Vol. 2, Progr. Math. 132,
(1996), 1--213.

\bibitem { } M. Gromov, Spaces and questions. GAFA 2000 (Tel Aviv, 1999).
 Geom. Funct. Anal. 2000, Special Volume, Part I, 118--161.

\bibitem { } M. Gromov and B. Lawson, Positive scalar curvature
and the Dirac operator on complete Riemannian manifolds.
 Inst. Hautes Études Sci. Publ. Math. No. 58 (1983), 83--196.

\bibitem { } E. Guentner, N. Higson and S. Weinberger, The
Novikov conjecture for linear groups. Preprint (2003).

\bibitem { } N. Higson and G. G. Kasparov, Operator K-theory for
groups which act properly and isometrically on Hilbert space,
Electronic Research Announcements, AMS 3 (1997), 131--141.

\bibitem {} N. Higson and  G. G. Kasparov, $E$-theory and $KK$-theory for
groups which act properly and isometrically on Hilbert space. Invent. Math.
144 (2001), no. 1, 23--74.

\bibitem {} N. Higson, G.  Kasparov and J. Trout, A Bott periodicity
theorem for infinite-dimensional Euclidean space. Adv. Math. 135
(1998), no. 1, 1--40.

\bibitem { } N. Higson and J. Roe, On the coarse Baum--Connes
conjecture, Novikov Conjectures, Index Theorems and Rigidity, Vol.
2, (S. Ferry, A. Ranicki and J. Rosenberg, editors), Cambridge
University Press, (1995), 227--254.

\bibitem { } N. Higson, J. Roe and G. Yu, A coarse
Mayer--Vietoris principle, Math. Proc. Camb. Philos. Soc. 114
(1993), 85--97.

\bibitem { } M. Hilsum and G. Skandalis, Invariance par homotopie de
la signature à coefficients dans un fibré presque plat.  J. Reine
Angew. Math. 423 (1992), 73--99.

\bibitem { } W. B. Johnson and N. L. Randrianarivony, $l_p$  ($p>2$)
does not coarsely embedd into a Hilbert space, preprint 2004.

\bibitem { } G. G. Kasparov, Equivariant KK-theory and the Novikov
conjecture, Invent. Math.  91 (1988), 147--201.

\bibitem { } G. G. Kasparov and G. Skandalis, Groups acting
properly on ``bolic'' spaces and the Novikov conjecture. Ann. of
Math. (2) 158 (2003), no. 1, 165--206.

\bibitem { } V. Lafforgue,  $K$-th\'eorie bivariante pour les alg\'ebres de
Banach et conjecture de Baum-Connes.  Invent. Math. 149 (2002),
no. 1, 1--95.

\bibitem { } M. Mendel and A. Naor, Metric Cotype, preprint, 2005.

\bibitem { } N. Ozawa, A note on non-amenability of $ B(l\sb p)$ for $p=1,2$.
Internat. J. Math. 15 (2004), no. 6, 557--565.

\bibitem { } J. Roe, Coarse cohomology and index theory for
complete Riemannian manifolds, Memoirs A.M.S., No. 104 (1993).

\bibitem { } J. Roe, Index Theory, Coarse Geometry, and Topology
of Manifolds, CBMS Regional Conf. Series in Math, Number 90, AMS
(1996).

\bibitem { } J. Rosenberg, $C^\ast$-algebras, positive scalar
curvature and the Novikov Conjecture, Publ I.H.E.S., No. 58
(1983), 197--212.

\bibitem { } C. Schochet, Topological methods for $C\sp{*} $-algebras. III.
Axiomatic homology. Pacific J. Math. 114 (1984), no. 2, 399--445.


\bibitem { } G. Yu, Coarse Baum--Connes conjecture, K-Theory 9
(1995), 199--221.

\bibitem { } G. Yu, Localization algebras and the coarse
Baum--Connes conjecture, K-Theory 11 (1997), 307--318.

\bibitem { } G. Yu, The Novikov Conjecture for groups with finite
asymptotic dimension, Ann. of Math.  Vol. 147, 2 (1998), 325--355.

\bibitem { } G. Yu, The coarse Baum-Connes conjecture for spaces which admit
 a uniform embedding into Hilbert space. Invent. Math. 139 (2000), no. 1, 201--240.

\bibitem { } G. Yu, Hyperbolic groups admit proper affine isometric actions on
$l^p$-spaces,  Geom. Funct. Anal. 15 (2005), no. 5.
\end{thebibliography}
\end{document}